\documentclass[review]{elsarticle}              

\usepackage{lineno,hyperref}
\modulolinenumbers[5]

\usepackage{graphicx,color}
\usepackage{amsmath,amsfonts,amssymb,amsthm}
\usepackage{ulem}

\definecolor{macouleur}{cmyk}{0.1,0.95,0,0.05}
\newcommand{\tran}{\mathsf{T}}

\newcommand{\bsy}{\boldsymbol{y}}
\newcommand{\bsmu}{\boldsymbol{\mu}}
\newcommand{\bsbeta}{\boldsymbol{\beta}}
\newcommand{\bstheta}{\boldsymbol{\theta}}
\newcommand{\bsSigma}{\boldsymbol{\Sigma}}
\newcommand{\bsX}{\boldsymbol{X}}

\newcommand{\bsU}{\boldsymbol{U}}
\newcommand{\dnorm}{\mathcal{N}}
\newcommand{\dustd}{\mathcal{U}}
\newcommand{\real}{\mathbb{R}}
\newcommand{\ult}{\underline{\tau}}

\journal{xxx}

\begin{document}

\begin{frontmatter}

\title{Shapley effects for sensitivity analysis with correlated inputs: comparisons with Sobol' indices, numerical estimation and applications}

\author{Bertrand Iooss}
\address{EDF Lab Chatou, 6 Quai Watier, 78401 Chatou, France}
\address{Institut de Math\'ematiques de Toulouse, 31062, Toulouse, France}

\author{Cl\'ementine Prieur}
\address{Universit\'e Grenoble Alpes, CNRS, LJK, F-38000 Grenoble, France}
\address{Inria project/team AIRSEA}

\begin{abstract}
The global sensitivity analysis of a numerical model aims to quantify, by means of sensitivity indices estimates, the contributions of each uncertain input variable to the model output uncertainty.
The so-called Sobol' indices, which are based on functional variance analysis, present a difficult interpretation in the presence of statistical dependence between inputs.
The Shapley effects were recently introduced to overcome this problem as they allocate the mutual contribution (due to correlation and interaction) of a group of inputs to each individual input within the group.
In this paper, using several new analytical results, we study the effects of linear correlation between some Gaussian input variables on Shapley effects, and compare these effects to classical first-order and total Sobol' indices.
This illustrates the interest, in terms of sensitivity analysis setting and interpretation, of the Shapley effects in the case of dependent inputs. 
For the practical issue of computationally demanding computer models, we show that the substitution of the original model by a metamodel (here, kriging) makes it possible to estimate these indices with precision at a reasonable computational cost.
\end{abstract}


\end{frontmatter}

\section{Introduction}

When constructing and using numerical models simulating physical phenomena, global sensitivity analysis (SA) methods are valuable tools \cite{ioolem15,weilu15,borpli16,ioosal17}. 
These methods allow one to determine which model input variables contribute the most to the variability of the model outputs, or on the contrary which are not important and possibly which variables interact with each other.
The standard quantitative methods compute the variance-based sensitivity measures  also called Sobol' indices.
In the simple framework of $d$ scalar inputs, denoted by $\bsX=(X_1,\ldots,X_d)^\tran \in \real^d$ and a single scalar output $Y \in \real$, the model response is
\begin{equation}\label{eq:model}
Y = f(\bsX) \;.
\end{equation}
In the case of independent inputs, the interpretation of the Sobol' indices is simple because the variance decomposition of $Y$ is unique \cite{sob93}. 
For instance, the first-order Sobol' index of $X_i$, denoted $S_i$, represents the amount of the output variance solely due to $X_i$. The second-order Sobol' index ($S_{ij}$) expresses the contribution of the interactions of the pairs of variables $X_i$ and $X_j$, and so on for the higher orders. 
As the sum of all Sobol' indices is equal to one, the indices are interpreted as proportions of explained variance.

However, in many applications, it is common that the input variables have a statistical dependence structure imposed by a probabilistic dependence function \cite{kurcoo06} (e.g., a copula function) or by physical constraints upon the input or the output space \cite{petioo10,kuckly17,lopbol17}.
As shown in previous studies and presented in the following section, estimating and interpreting Sobol' indices with correlated inputs is not trivial \cite{saltar02,davwah08}.

{\color{black}The Shapley value, introduced in \cite{sha53}, is a solution concept in cooperative game theory. It has proved, however, to be a powerful tool in modelling
some economic problems. Indeed, since it can be interpreted
in terms of {\it marginal worth}, it is closely related to traditional
economic ideas. 
In \cite{owe14}, Art Owen has proposed to use that concept to measure variable importance in SA of model outputs. In the framework of dependent input variables, it allows namely to bypass the intricate issue of variance decomposition \cite{sonnel16,owepri17}. 
More precisely, the Shapley value is a solution concept of fairly distributing both gains and costs to several actors working in coalition. The Shapley value applies primarily in situations when the contributions of each actor are unequal. The Shapley value ensures each actor gains as much or more as they would have from acting independently. Now, if the actors are identified with a set of random model inputs and the value assigned to each coalition is identified to the explanatory power of the subset of model inputs composing the coalition, then the Shapley value can be interpreted as an importance measure of model inputs.

The two main properties and advantages of the Shapley values in the framework of variable importance measures are that they cannot be negative and they sum-up to the total output variance.
The equitable principle driving the allocation rule states that an interaction effect is equally apportioned to each input involved in the interaction. Indeed, if the model inputs are independent, the Shapley value associated to input $i$, $i \in \{1, \ldots d\}$, is given by 
\begin{equation}\label{share}
Sh_i=\sum_{u \subseteq \{1, \ldots , d\}, i \in u}\displaystyle \frac{S_u}{|u|}\;,
\end{equation}
where $S_u$ denotes the usual Sobol' index associated to $u$ and defined as $\displaystyle S_u=\frac{\sum_{v \subseteq u}(-1)^{|u|-|v|}\mbox{Var}(\mathbb{E}[Y|{\bsX}_v])}{\mbox{Var}(Y)}$, and $|u|$ is the cardinality of $u$

From the conceptual point of view, the major issue then remains to understand the effect of the dependence between inputs on the variance-based Shapley values.
Analytical computations for
several test-cases have been presented in \cite{owepri17}, e.g., a general formula was provided in the framework of Gaussian inputs $\bsX$ for a linear model $f(\bsX) = \beta_0+\bsbeta^\tran \bsX$ (with $\bsbeta \in \real^d$).
Moreover, the study in \cite{owepri17} has highlighted several properties of these indices (e.g., for $d  \geq 2$, if there is a bijection between any two of the $X_j$, $j=1, \ldots ,d$, then those two variables have the same Shapley value). 
However, from these analytical results, presented in Section \ref{sec:lingausgeneral}, it still seems  difficult to understand the effects of the input correlation structure onto the Shapley values.
Therefore, in this paper, we provide a thorough investigation of several particular cases,  simple enough to provide some guidance to their interpretation.}

For the sake of practical applications, \cite{sonnel16} has proposed two estimation algorithms of the Shapley effects (that we define as the normalized variance-based Shapley values), and illustrated them in two application cases.
As for the Sobol' indices \cite{ioolem15,pritar17}, one important issue in practice is the numerical cost in terms of number of model evaluations required to estimate the Shapley effects.
To ease the computational burden, a classical solution is to use a metamodel which is a mathematical approximation of the numerical model (\ref{eq:model}) from an initial and limited set of runs \cite{fanli06,ioovan06}.
The metamodel solution is a usual engineering practice for estimating sensitivity indices \cite{legmar17}.
In this paper, we investigate the practical use of metamodeling for estimating the Shapley effects.

In the following section, we synthesize the previous works and recall the general mathematical formulation of Sobol' indices and Shapley effects when the inputs are dependent.
We also discuss the SA settings \cite{saltar02} that can be addressed with the Shapley effects.
In Section \ref{sec:theory}, we develop the analytical formulas that one can get in several particular cases: linear models with Gaussian inputs and block-additive structure in two and three dimensions. 
In particular, we focus on inequalities that can be established between Sobol' indices and Shapley effects.
Metamodel-based algorithms for estimating  Shapley effects are considered in Section \ref{sec:estimation} which also presents an industrial application.
The conclusion summarizes the findings and contributions of this work.

\section{General formulation of sensitivity indices}\label{sec:general}

\subsection{Sobol' indices}\label{sob:ind}

Starting from the model (\ref{eq:model}) $Y=f(\bsX)$, the Sobol' index associated with a set of inputs indexed by $u$ (${u} \subseteq \{1, \ldots , d\}$) is defined by:
\begin{equation} \label{indclosed}
S_u^{\textup{clo}} = \mbox{Var}(\mathbb{E}[Y|\bsX_u])/\mbox{Var}(Y) \;.
\end{equation}
Sobol' indices have been introduced in \cite{sob93}. Indices defined by \eqref{indclosed} are referred as closed Sobol' indices in the literature (see, e.g., \cite{pritar17}) and are always comprised between $0$ and $1$.
For the mathematical developments of the following sections, we denote the numerator of $S_u^{\textup{clo}}$ as
\begin{equation}\label{eq:tau}
\underline{\tau}_{{u}}^2 = \mbox{Var}(\mathbb{E}[Y|\bsX_{{u}}])\;.
\end{equation}

In addition to the Sobol' indices in Eq. (\ref{indclosed}), total sensitivity indices have also been defined in order to express the ``total'' sensitivity of the variance of $Y$ to an input variable $X_i$ \cite{homsal96}:
\begin{equation}\label{indtotal2}
S_{T_i} = \frac{\mathbb{E}_{\bsX_{-i}}(\mbox{Var}_{X_i}[Y|\bsX_{-i}])}{\mbox{Var}(Y)},
\end{equation}
{\color{black}where $\bsX_{-i}$ is the $(d-1)$-dimensional vector $(X_1, \ldots, X_{i-1},X_{i+1}, \ldots , X_d)$. Note that $\mathbb{E}_{\bsX_{u}}$ ({\it resp.} $\mbox{Var}_{\bsX_u}$) means the expectation ({\it resp.} the variance) with respect to the random vector$\bsX_u$.}  

\subsection{Sensitivity indices with dependent inputs}

Recall that this paper is not restricted to independent inputs, so that the knowledge of first-order and total Sobol' indices does not give complete information on the way an input $X_i$ influences the output $Y$. 
Many propositions appear in the literature that are not always easy to interpret.
One strategy proposed by \cite{jaclav06} is to evaluate the Sobol' indices of subsets of inputs which are correlated within the subsets but not correlated outside. 
However, this approach is not satisfactory because one may need to compute the Sobol' indices of the individual variables.

In \cite{xuger08} the authors proposed to decompose each partial variance $V_i$ due to input $X_i$, and defined as $\mbox{Var}(\mathbb{E}[Y|X_i])$,  into partial variances ($V_i^{\textup{U}}$) due to the uncorrelated variations of input $X_i$ and partial variance ($V_i^{\textup{C}}$) due to the correlated variations of $X_i$ with all other inputs $X_j$, $j\neq i$.
 Such an approach allows to exhibit inputs that have an impact on the output only through their strong correlation with other inputs. 
However, their approach only applies to linear model outputs with linearly dependent inputs. 
In the same spirit, the authors in \cite{martar12} proposed a strategy which supports non-linear models and non-linear dependencies. 
Their methodology is decomposed in two steps: a first step of decorrelation of the inputs and then a second step based on the concept of High Dimensional Model Representation (HDMR).  
HDMR (see, e.g., \cite{lira08}) relies on a hierarchy of component functions of increasing dimensions. The second step in \cite{martar12} thus consists in performing the HDMR on the decorrelated inputs.
At the same time, the authors in \cite{kuctar12} proposed a non-parametric procedure to estimate first-order and total indices in the presence of dependencies, not necessarily of linear type. 
Their methodology requires knowledge of the conditional probability density functions and the ability to draw random samples from those. 
Later, in \cite{martar15}, the authors established a link between the approaches in \cite{kuctar12} and \cite{martar12}, allowing the distinction between the independent contributions of inputs to the response variance and their mutual dependent contributions, via the estimation of four sensitivity indices for each input, namely full and independent first-order indices, and full and independent total indices (see the following section). 

A different approach was introduced in \cite{lirab10}. It is once again based on the HDMR. The component functions are approximated by expansions in terms of some suitable basis functions (e.g., polynomials \cite{sudcan13}, splines \ldots). 
This meta-modeling approach allows a covariance decomposition of the response variance.
In \cite{zhalu15}, the output variance is decomposed into orders of partial variance contributions, while the second order and higher orders of partial variance contributions are decomposed into uncorrelated interaction contributions and correlated contributions.

It is worth noting that none of these works has given an univocal definition of the functional ANOVA for correlated inputs as the one provided by the Hoeffding-Sobol' decomposition \cite{hoe48,sob93} when inputs are independent. 
A generalization of the Hoeffding-Sobol' decomposition was proposed in Stone \cite{sto94} (see also Hooker \cite{hoo07}). 
Then, in \cite{chagam12}, the authors defined a new variable importance measure based on the decomposition of Stone. 
However, such an important measure suffers from two conceptual problems underlined in \cite{owepri17}: the sensitivity indices can be negative and the approach places strong restrictions on the joint probability distribution of the inputs.

As a different approach, \cite{bor07,borcas11} initiated the construction of novel generalized moment-free sensitivity indices, called $\delta$-importance measures. 
Based on some geometrical considerations,
these indices measure the shift area between the outcome probability density function and this
same density conditioned to an input. Thanks to the properties of these
new indices, a methodology is available to obtain them analytically through test
cases.
In \cite{lopbol17} an application of these sensitivity measures on a gas transmission model with dependent inputs is proposed. 
We note that the authors in \cite{zholu14} have proposed a methodology to evaluate the $\delta$-importance measures by employing the decorrelation procedure described in \cite{martar12}.

\subsection{Full and independent Sobol' indices}

In \cite{martar15}, the authors propose a strategy based on the estimation of four sensitivity indices per input, and described hereafter. 
{\color{black}We assume that $\bsX=\left(X_1, \ldots , X_d\right)$ is a set of continuous dependent random variables, with joint probability density function $p(\bsX)=p(X_1)p(X_2|X_1)p(X_3|X_1,X_2)\ldots p(X_d|\bsX_{-d})$. The Rosenblatt transform introduced in \cite{ros52} is the most general transformation to map the $d$-dimensional input vector $\bsX$ onto a $d$-dimensional vector with independent components.
Let $F_i(x_i|\bsX_u)$ denote the cumulative distribution function of $X_i$ conditionally to $\bsX_u$, with $u \subseteq \{1, \ldots, i-1,i+1, \ldots,d\}$. 
The Rosenblatt transform of vector $\bsX$ is defined as $T_{\textup{Ros}, 1}(\bsX)=\bsU^1=(U_1^1, \ldots , U_d^1)$ with
$U_1^1=F_1(X_1)$, $U_2^1=F_2(X_2|X_1)$, $U_3^1=F_3(X_3|X_1,X_2)$, $\ldots$ , $U_d^1=F_d(X_d|\bsX_{-d})$. 
In case of dependent variables, the Rosenblatt transform is not unique and there are $n!$ possibilities depending on how the random variables are ordered in the vector $\bsX$. 
We denote by $T_{\textup{Ros}, i}(\bsX)=\bsU^i$ the Rosenblatt transform of the ordered vector $(X_i,X_{i+1},\ldots,X_d,X_1, \ldots ,X_{i-1})$. 
The transforms $T_{\textup{Ros},i}$, $i=1, \ldots, d$, are bijective and for any $i=1, \ldots , d$ we can write $Y=f(X_1, \ldots , X_d)=g_{i}(\bsU^{i})$ with $g_i=f \circ T_{\textup{Ros},i}^{-1}$. 
We now define, as in \cite{martar15},

\begin{eqnarray}
\displaystyle S_{(i)} & = & \frac{\mbox{Var}(\mathbb{E}[g_{i}(\bsU^{i})|U^{i}_1])}{\mbox{Var}(g_i(\bsU^i))}= \frac{\mbox{Var}(\mathbb{E}[Y|X_i])}{\mbox{Var}(Y)}=S_i \, ,\label{eq:1}\\
\displaystyle S_{T_{(i)}} & = & \frac{\mathbb{E}(\mbox{Var}[g_{i}(\bsU^{i})|\bsU^{i}_{-1}])}{\mbox{Var}[g_i(\bsU^i)]} \, ,\label{eq:2}\\
\displaystyle S_{(i)}^{\textup{ind}} & = & \frac{\mbox{Var}(\mathbb{E}[g_{i+1}(\bsU^{i+1})| U^{i+1}_d])}{\mbox{Var}[g_{i+1}(\bsU^{i+1})]} \, ,\label{eq:3}\\
\displaystyle S_{T_{(i)}}^{\textup{ind}} & = & \frac{\mathbb{E}(\mbox{Var}[g_{i+1}(\bsU^{i+1})|\bsU^{i+1}_{-d}])}{\mbox{Var}[g_{i+1}(\bsU^{i+1})]}= \frac{\mathbb{E}(\mbox{Var}[Y|\bsX_{-i}])}{\mbox{Var}(Y)}=S_{T_i} \, .\label{eq:4}
\end{eqnarray}}

The indices $S_{(i)}$ and $S_{T_{(i)}}$ include the effects of the dependence of $X_i$ with other inputs, and are referred to as full sensitivity indices in \cite{martar12}. 
The indices $S_{(i)}^{\textup{ind}}$ and $S_{T_{(i)}}^{\textup{ind}}$ measure the effect of an input $X_i$, that is not due to its dependence with other variables $\bsX_{-i}$. 
Such indices have also been introduced as uncorrelated effects in \cite{martar12} and further discussed in \cite{martar15} which refers to them as the independent Sobol' indices.
In \cite{martar15}, the authors propose to estimate the four indices $S_{(i)}$, $S_{T_{(i)}}$, $S_{(i)}^{\textup{ind}}$ and $S_{T_{(i)}}^{\textup{ind}}$ for the full set of inputs ($i=1, \ldots , d$). 
Note that $S_{(i)}^{\textup{ind}} \leq S_{T_{(i)}}^{\textup{ind}}=S_{T_{i}}$ and that $S_i = S_{(i)} \leq S_{T_{(i)}}$, but other inequalities are not known.
\cite{martar15} has proposed two sampling strategies for dependent continuous inputs. The first one is based on the Rosenblatt transform \cite{ros52}. The second one is a simpler method that estimates the sensitivity indices without requiring the knowledge of conditional probability density functions, and which can be applied in case the inputs dependence structure is defined by a rank correlation matrix (see, e.g., \cite{imacon82}).

In the following, we focus our attention on the full first-order Sobol' indices $S_{(i)}$ and the independent total Sobol' indices $S_{T_{(i)}}^{\textup{ind}}$.
Indeed, these are the indices which are used in the definition of the Shapley effects (see Section \ref{sec:shapley}).
Moreover, the full first-order indices $S_{(i)}$ coincide with the associated classical first-order Sobol' indices $S_{i}$, and the independent total indices $S_{T_{(i)}}^{\textup{ind}}$ coincide with the associated classical total Sobol' indices $S_{T_{i}}$. 

\subsection{Shapley effects}\label{sec:shapley}

{\color{black}As mentioned in the introduction, the Shapley value, introduced in \cite{sha53}, is a solution concept in cooperative game theory, which has been adapted to measure variable importance in \cite{owe14}. 
The Shapley value ensures that each input variable has as much sensitivity as it would when acting independently.
The value assigned to each coalition (set of inputs) is identified to its explanatory power: more precisely, we consider in the following that the value assigned to coalition $\bsX_u$ is \begin{equation}\label{eq:costfct}
c(u)=S_u^{\textup{clo}}=\mbox{Var}(\mathbb{E}[Y|\bsX_u])/\mbox{Var}(Y) \;\cdot
\end{equation}
The Shapley values corresponding to the cost function $c(\cdot)$ defined above are  known as Shapley effects in the literature. The equitable principle driving the Shapley sharing rule leads then to the definition of Shapley effects $Sh_i$ as
\begin{equation}\label{indshapley}
Sh_i = \sum_{u \subseteq -\{i\}} \frac{(d-|u|-1)! |u|!}{d!} \left[ c(u \cup \{i\}) - c(u) \right] \;,
\end{equation}
with $-\{i\}$ the set of indices $\{1, \ldots, d\}$ not containing $i$.
The Shapley effects then consist in importance measures, which can be used for SA of model output.}
The authors in \cite{sonnel16} also prove that it is equivalent to define $c(\cdot)$ as $
\mathbb{E}[\mbox{Var}(Y|\bsX_{-u})]/\mbox{Var}(Y)$ or as $\mbox{Var}(\mathbb{E}[Y|\bsX_u])/\mbox{Var}(Y)$.
Note that, in \cite{owe14} and \cite{sonnel16}, the cost function is not normalized by the variance of $Y$ while, in  the present paper, we consider its normalized version.
 
Recall that the Shapley effects result from an equitable sharing of the model output variance between the input variables.
In the framework of independent input variables, it means that the Shapley effect $Sh_i$ associated to input $X_i$ shares the part of variance due to the interaction of a subset of inputs $\bsX_u$, $u \subseteq \{1, \ldots  d\}$, $i \in u$ with each individual input within that subset (see Equation \eqref{share} in the introduction).
In the framework of dependent input variables, the Shapley effect associated to input $X_i$ ($i \in \{1, \ldots , d\}$) takes into account both interactions and correlations of $X_i$ with $X_j$, $1 \leq j \leq d$, $j \neq i$. 
The consequence is that Shapley effects are non negative and sum-up to one, allowing an easy interpretation for ranking input variables.

Formulas (\ref{eq:costfct}-\ref{indshapley}) show that the Shapley effect of an input is a by-product of its Sobol' indices. 
Thus, if one can compute the complete set of full and independent indices (Eqs. (\ref{eq:1}-\ref{eq:4})), we can compute the Shapley effect of each input. 
Two estimation algorithms, based on consistent estimators of the Shapley effects, have been proposed in \cite{sonnel16}: The ``Exact permutation'' algorithm traverses all possible permutations between the inputs and the ``Random permutation'' algorithm consists of randomly sampling some permutations of the inputs (see \ref{sec:direct} for more details).
From the exact permutation algorithm, we can extract a consistent estimator of any Sobol' index. Concerning the random permutation algorithm, the sample size $N_i$ related to the inner loop (conditional variance estimation), and the one $N_o$ related to the outer loop (expectation estimation) are fixed to $N_I=3$ and $N_o=1$ respectively (following the recommendation of \cite{sonnel16}, see also \ref{sec:direct}).
Thus it is not possible to extract from that algorithm accurate estimates of Sobol' indices. However, this last algorithm is consistent for the estimation of Shapley effects and is particularly adapted in the case of high-dimensional input spaces.

In case the input variables are independent, first-order ({\it resp.} total) Sobol' indices provide effectively computable lower- ({\it resp.} upper-) bounds for the Shapley effects. 
In the following, as in \cite{sonnel16}, we will show that these bounds do not hold anymore in case the input factors present some dependencies. 

\subsection{SA settings}\label{sec:SAsettings}

{\color{black} In this section, a ``factor'' is the term, coming from the field of experimental designs, referring to an input variable.}
\cite{saltar02} and \cite{saltar04} have defined several objectives, called SA settings, that sensitivity indices can address.
These SA settings aims at clarifying the objectives of the analysis.
They are listed in \cite{saltar04} as follows:
\begin{itemize}
\item Factors Prioritization (FP) Setting, to know on which inputs the reduction of uncertainty leads to the largest reduction of the output uncertainty;
\item Factors Fixing (FF) Setting, to determine which inputs can be fixed at given values without any loss of information in the model output;
\item Variance Cutting (VC) Setting, to know which inputs have to be fixed to obtain a target value on the output variance;
\item Factors mapping (FM) Setting, to determine which inputs are most responsible for producing values of the output in a specific region of interest.
\end{itemize}

In the case of independent inputs, the Sobol' indices directly address the FP, FF and VC Settings (see \cite{saltar04}).
{\color{black}In the framework of dependent inputs, the FP and FF settings have to be clarified. We consider in the present work that the variance of one variable cannot be changed independently of the remaining variables, due to the correlation structure of the vector $\bsX$.
A direct consequence is that the classical ANOVA-Sobol' decomposition does not hold true anymore and thus the FP and VC settings cannot be directly obtained with Sobol' indices (see  \cite{saltar02}).
The FF Setting is also more difficult to address in the presence of dependencies. More precisely, we mean by FF that one or more factors is fixed while the remaining factors are sampled from the conditional distribution with respect to the values of the fixed variables. Thus, fixing one or more of the input factors has an impact on all the input factors that are correlated with them. Then the expectation over the different values of the fixed variables should be considered.}
{\color{black} However, as explained in \cite{martar15}, independent and full total Sobol' indices can help understanding the FF Setting in the dependent framework. Indeed, if a factor has both indices $S_{T_{(i)}}^{\textup{ind}}$ and $S_{T_{(i)}}$ which are null, then we can deduce that the model output $Y$ can be expressed only as a function of $\bsX_{-i}$ or $\bsX_{-i}|X_i$ (see \cite{beneli18} for more details).}

It is of interest now to give some hints about how modelers can use the Shapley effects to address some SA settings.
The VC Setting is not achieved, in the independent and dependent inputs cases, because the Shapley effect of an input contains some effects due to other inputs.
In the case of independent inputs, each Shapley effect $Sh_i$ is bounded by the corresponding first-order $S_i=S_{(i)}$ and total indices $S_{T_i}=S_{T_{(i)}}^{\textup{ind}}$:
\begin{equation}\label{sand}
S_i \leq Sh_i \leq S_{T_{i}}\;.
\end{equation}
In addition to the individual effect of the variable $X_i$, the Shapley effects take into account the effects of interactions by distributing them equally in the index of each input that plays in the interaction \cite{owe14} {\color{black}(see Equation \eqref{share} in the introduction)}.
Therefore, the FF Setting is achieved by using the Shapley effects. However, the FP Setting is not precisely achieved because we cannot distinguish the contributions of the main and interaction effects in a Shapley effect for the case of independent inputs.

In the case of dependent inputs, Eq. \eqref{sand} does not hold true anymore. However, due to the equitable principle on which the allocation rule is based, a Shapley effect close to zero means that the input has no significant contribution to the variance of the output, neither by its interactions nor by its dependencies with other inputs. Therefore, the FF Setting can be addressed with $Sh_i$ for the case of non-independent inputs, which cannot be directly done with either of Sobol' indices.

\section{Relations and inequalities between Sobol' indices and Shapley effects}\label{sec:theory}

As said before, in the case of dependent inputs, no relation such as the one in Eq. \eqref{sand} can be directly deduced.
The goal of this section is to study particular cases where analytical deductions can be made.
The linear model is the first model to be studied in sensitivity analysis because of its simplicity, its interpretability, and the tractable analytical results provided for Gaussian inputs. 
Considering two and three inputs allows comparing our new results to previous studies on Sobol' indices.
The Ishigami function is the most used analytical test function in sensitivity analysis.
It has been considered in numerous papers dealing with Sobol' indices with dependent inputs.

We focus the analysis on the full first-order indices (corresponding to classical first-order Sobol' indices) $S_i=S_{(i)}$ (called ``First-order Sobol'' in the figures) and the independent total  indices (corresponding to classical total Sobol' indices) $S_{T_i}=S_{T_{(i)}}^{\textup{ind}}$ (called ``Ind total Sobol'' in the figures).
Indeed these are the indices mainly studied and discussed in the previous works on this subject, in particular in \cite{kuctar12}.
The numerical tests of \cite{kuctar12} have inspired the ones proposed in this section.
Moreover, these indices are easily and directly provided by the Shapley effects estimation algorithms (\cite{sonnel16}, see Section \ref{sec:estimation}), during the first and last iterations of the algorithm.
Finally, the estimation of the independent first-order Sobol' $S_{(i)}$ and of the full total Sobol' indices $S_{T_{(i)}}$ is based on a rather cumbersome process, the use of $d$ Rosenblatt transforms (see more details in Section \ref{sob:ind}).
Comparisons with these complementary indices are made in \cite{beneli18}.

\subsection{Linear model with Gaussian input variables}\label{sec:lingausgeneral}

Let us consider 
\begin{equation}\label{eq:linearmodel}
Y=\beta_0+\bsbeta^\tran\bsX
\end{equation}
with $\bsX\sim \dnorm(\bsmu,\bsSigma)$ and $\bsSigma \in \real^{d \times d}$ a positive-definite matrix.
We have $\sigma^2=\mbox{Var}(Y)=\bsbeta_{\{1,\ldots,d\}}^\tran\bsSigma_{\{1, \ldots , d\},\{1, \ldots , d\}}\bsbeta_{\{1,\ldots,d\}}$.
Note that the subscripts are added on $\bsbeta$ (resp. $\bsSigma$) in order to identify which components are contained in the vector (resp. matrix).

We get from \cite{owepri17}:
\begin{equation}\label{eq:shgaus}
Sh_i = \frac1d\sum_{u\subseteq-i}{d-1\choose |u|}^{-1} \frac{ \mbox{Cov}\bigl(X_i,\bsX_{-u}^\tran\bsbeta_{-u}\mid\bsX_u\bigr)^2} {\sigma^2 \, \mbox{Var}(X_i\mid\bsX_u)} \, .
\end{equation}
Recall now the following classical formula:
{\color{black}\begin{equation}\label{varcond}
\mbox{Var}(\bsX \mid\bsX_{-j})=\bsSigma_{\{1, \ldots , d\},\{1, \ldots , d \}}-\bsSigma_{\{1, \ldots , d\},-j}\bsSigma_{-j,-j}^{-1}\bsSigma_{-j,\{1, \ldots , d\}}\, .
\end{equation}}
From \eqref{varcond} and according to the law of total variance, we easily obtain the Sobol' indices:
{\color{black}\begin{align}\displaystyle
S_j
& \displaystyle = \frac{\mbox{Var}(\mathbb{E}[\beta_0+\bsbeta^\tran\bsX|X_j])}{\sigma^2}=1-\frac{\mathbb{E}(\mbox{Var}[\bsbeta^\tran\bsX|X_j])}{\sigma^2}
\nonumber \\ 
& \displaystyle = \frac{\bsbeta_{\{1,\ldots,d\}}^\tran\bsSigma_{\{1, \ldots , d\},j}\bsSigma_{j,j}^{-1}\bsSigma_{j,\{1, \ldots , d\}}\bsbeta_{\{1, \ldots , d\}}}{\sigma^2} \label{o1}\, ,
\end{align}}

{\color{black}\begin{align}
\displaystyle
S_{T_{j}} 
&\displaystyle = \frac{\mathbb{E}(\mbox{Var}[\beta_0+\bsbeta^\tran\bsX|\bsX_{-j}])}{\sigma^2}=\frac{\mathbb{E}(\mbox{Var}[\bsbeta^\tran\bsX|\bsX_{-j}])}{\sigma^2} \nonumber \\
&\displaystyle = \frac{\bsbeta_{\{1,\ldots,d\}}^\tran\left(\bsSigma_{\{1, \ldots , d\}, \{1, \ldots , d \}}-\bsSigma_{\{1, \ldots , d\},-j}\bsSigma_{-j,-j}^{-1}\bsSigma_{-j,\{1, \ldots , d\}}\right)\bsbeta_{\{1, \ldots , d\}}}{\sigma^2} \label{ot}\, \cdot
\end{align}}

Note that $\beta_0$ and $\bsmu$ do not play any role as translation parameters in variance-based sensitivity analysis.
However, no general conclusion can be drawn from Eqs. (\ref{eq:shgaus}-\ref{ot}). 
Therefore, particular cases are studied in the three following sections.

\subsection{Linear model with two Gaussian inputs}\label{sec:lin2gau}

Let us consider the case $d=2$ with 
$$\bsmu=\begin{pmatrix} \mu_1\\
\mu_2
\end{pmatrix},
 \; \bsbeta=\begin{pmatrix} \beta_1\\
\beta_2
\end{pmatrix} \mbox{ and }\bsSigma=
\begin{pmatrix} 
\sigma_1^2 & \rho \sigma_1\sigma_2 \\
\rho \sigma_1\sigma_2  & \sigma_2^2
\end{pmatrix}, \; -1 \leq \rho \leq 1 \, , \sigma_1>0 \, , \sigma_2>0 \, .$$
We have $\displaystyle \sigma^2= \mbox{Var} (Y)=\beta_1^2 \sigma_1^2+2 \rho \beta_1 \beta_2  \sigma_1 \sigma_2+ \beta_2^2 \sigma_2^2$.
From Eq. (\ref{eq:tau}), $\underline{\tau}_{{u}}^2=\mbox{Var}(\mathbb{E}[Y|\bsX_{{u}}])$ (${u} \subseteq \{1,2\}$) and we obtain $\underline{\tau}^2_{\emptyset}=0$, $\underline{\tau}^2_1=(\beta_1\sigma_1+\rho \beta_2 \sigma_2)^2$, $\underline{\tau}^2_2=(\beta_2 \sigma_2+\rho \beta_1 \sigma_1)^2$ and $\underline{\tau}^2_{12}=\sigma^2$.
The definition of the Shapley effect (Eq. (\ref{indshapley})) gives ($j=1,2$)
\begin{align}
\sigma^2 Sh_j 
&= \frac1d\sum_{u\subseteq -\{j\}}
{d-1\choose |u|}^{-1}
(\ult^2_{u+\{j\}}-\ult^2_u),
\end{align} 
from which we get
\begin{equation}\label{shap}
\begin{array}{rcl}
\sigma^2 Sh_1 & = & \displaystyle \beta_1^2 \sigma_1^2 (1- \frac{\rho^2}{2})+\rho \beta_1 \beta_2 \sigma_1 \sigma_2+\beta_2^2 \sigma_2^2 \frac{\rho^2}{2} \, ,\\ 
\sigma^2 Sh_2 & = & \displaystyle \beta_2^2 \sigma_2^2 (1- \frac{\rho^2}{2})+\rho \beta_1 \beta_2 \sigma_1 \sigma_2+ \beta_1^2 \sigma_1^2 \frac{\rho^2}{2} \, .
\end{array}
\end{equation}

From Eq. (\ref{o1}-\ref{ot}), we have
\begin{equation}\label{sfull}
\begin{array}{rcl}
\sigma^2 \, S_1 & = & \displaystyle \beta_1^2  \sigma_1^2 + 2 \rho \beta_1 \beta_2 \sigma_1 \sigma_2 +  \rho^2 \beta_2^2 \sigma_2^2 \, , \\
\sigma^2 \, S_2 & = & \displaystyle \beta_2^2 \sigma_2^2 + 2 \rho \beta_1 \beta_2 \sigma_1 \sigma_2 + \rho^2 \beta_1^2 \sigma_1^2\, .
\end{array}
\end{equation}
and
\begin{equation}\label{tind}
\begin{array}{rcl}
\sigma^2 S_{T_{1}} & = & \displaystyle \beta_1^2 \sigma_1^2 (1- \rho^2) \, ,  \\
\sigma^2 S_{T_{2}}& = & \displaystyle \beta_2^2 \sigma_2^2 (1- \rho^2)  \, .
\end{array}
\end{equation}

From Equations (\ref{shap}), (\ref{sfull}) and (\ref{tind}) we can infer that the following four assertions are equivalent
\begin{equation}\label{eq:sandwich1}
\begin{array}{l}
Sh_j \leq S_{T_{j}} \, , \\
S_j \leq Sh_j \, , \\
\displaystyle \rho \, \left( \rho \frac{\beta_1^2 \sigma_1^2+ \beta_2^2 \sigma_2^2}{2}+ \beta_1 \beta_2 \sigma_1\sigma_2\right) \leq 0  \, , \\
\displaystyle |\rho| \leq  \frac{2 |\beta_1 \beta_2| \sigma_1 \sigma_2}{\beta_1^2 \sigma_1^2+ \beta_2^2 \sigma_2^2} \, .
\end{array}
\end{equation}
The equality of the first three assertions is obtained in the absence of correlation ($\rho=0$).
In that case, the Shapley effects are equal to the first-order and total Sobol' indices.
In the presence of correlation, the Shapley effects lie between the full first-order indices and the independent total indices: with either $S_j \leq Sh_j \leq S_{T_{j}}$ or $S_{T_{j}} \leq Sh_j \leq S_j$.
We call this the ``sandwich effect''.
We remark that the effects of correlations on the independent total indices ({\it e.g.} $-\rho^2 \beta_1^2 \sigma_1^2$ for $X_1$) and on the full first-order indices ({\it e.g.} $2\rho \beta_1 \beta_2 \sigma_1 \sigma_2 +  \rho^2 \beta_2^2 \sigma_2^2$ for $X_1$) are allocated half to the Shapley effect, in addition to the elementary effect ({\it e.g.} $\beta_1^2 \sigma_1^2$ for $X_1$).

\begin{figure}[!ht]
  \centering
	\includegraphics[width=0.49\textwidth]{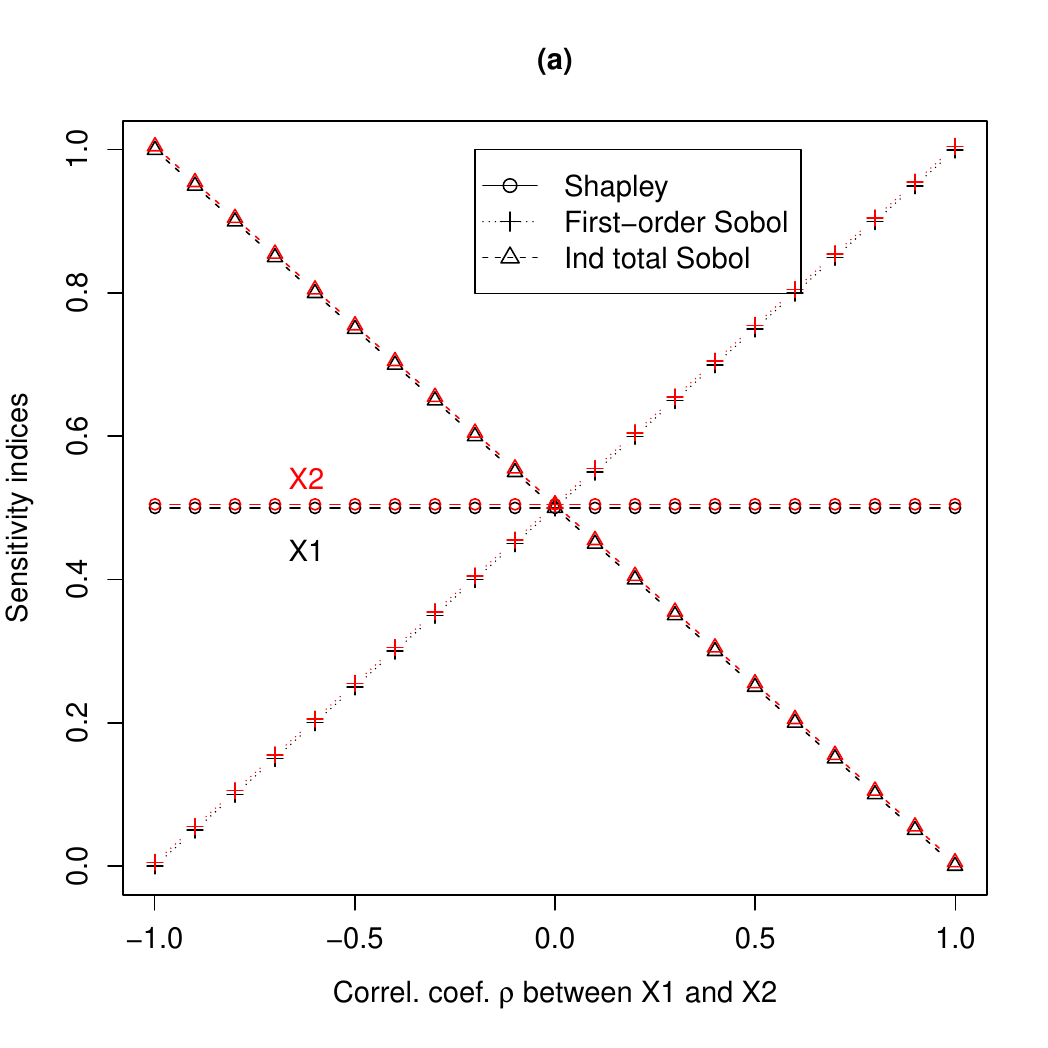}
	\includegraphics[width=0.49\textwidth]{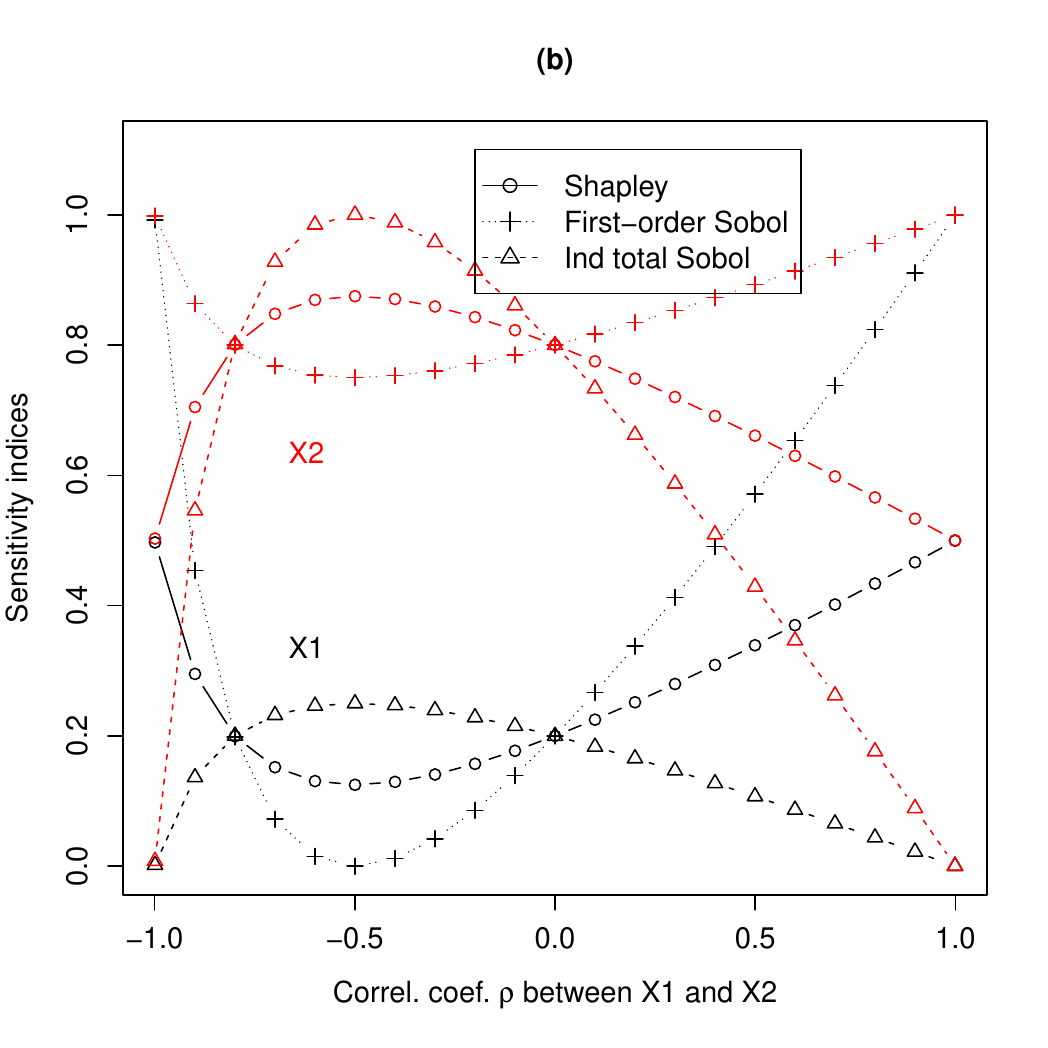}
		
\vspace{-0.4cm}
\caption{Sensitivity indices on the linear model ($\beta_1=1$, $\beta_2=1$) with two Gaussian inputs. (a): $(\sigma_1,\sigma_2) = (1,1)$. (b): $(\sigma_1,\sigma_2) = (1,2)$.}
  \label{fig:modlin2d}
\end{figure}

These results are also illustrated in Figure \ref{fig:modlin2d}.
In Figure \ref{fig:modlin2d} (a), as the standard deviations of each variable are equal, the different sensitivity indices are superimposed and the Shapley effects are constant.
In Figure \ref{fig:modlin2d} (b), because $X_2$ is more uncertain than $X_1$, its sensitivity indices are logically larger than those of $X_1$.
The effect of the dependence between the inputs is clearly shared on each input variable.
The dependence between the two inputs lead to a rebalancing of their corresponding Shapley effects, while a full Sobol' index of an input comprises the effect of another input on which it is dependent.
We also see on Figure \ref{fig:modlin2d} (b) that the Shapley effects of two perfectly correlated variables are equal.
Finally, the ``sandwich effect'' is respected for each input: From Eq (\ref{eq:sandwich1}), we can prove that $S_j \leq Sh_j \leq S_{T_{j}}$ when $\rho \in [-0.8;0]$ and that $S_j > Sh_j > S_{T_{j}}$ elsewhere.

\subsection{Correlated input non included in the model}\label{sec:lingaus2particular}

Consider the model $Y=f(X_1,X_2)=X_1$ with $(X_1,X_2)$ two dependent standard Gaussian variables with a correlation coefficient $\rho$.
It corresponds to the case $\beta_1=1$, $\beta_2=0$, $\mu_1=\mu_2=0$, $\sigma_1^2=\sigma_2^2=1$ in the model introduced in Section \ref{sec:lin2gau}.
The Shapley effects are
\begin{equation}\label{eq:shapley_modX1}
Sh_1 = 1 - \frac{\rho^2}{2} \mbox{ and } Sh_2 = \frac{\rho^2}{2}\;.
\end{equation}
Eq. (\ref{eq:shapley_modX1}) leads to the important remark that an input not involved in the numerical model can have a non-zero effect if it is correlated with an influential input of the model.
If the two inputs are perfectly correlated, their Shapley effects are equal.
This example also illustrates the FF setting that can be achieved with the Shapley effects: if $\rho$ is close to zero, $Sh_2$ is small and $X_2$ can be fixed without changing the output variance.

For the Sobol' indices, we have 
\begin{equation}\label{eq:sobol_modX1}
S_1 = 1 \;,\; S_{T_1} = 1 - \rho^2 \mbox{ and } S_2 = \rho^2\;,\; S_{T_2} = 0 \;,
\end{equation}
which indicates that $X_2$ is only important because of its correlation with $X_1$ (FP setting) and that by only accounting for the uncertainty in $X_1$, one should be able to evaluate the uncertainty of $Y$ accurately. {\color{black}More generally, if a black box model has $d$ inputs, and if $S_{T_i}=0 $, it means that the model output can be written as a measurable function of $(X_1, \ldots , X_{i-1},X_{i+1}, \ldots , X_d)$ only. Thus, if then $S_i>0$, it means that $X_i$ is correlated to $(X_1, \ldots , X_{i-1},X_{i+1}, \ldots , X_d)$.}

\subsection{Linear model with three Gaussian inputs}\label{sec:lin3gau}

We consider the linear model (\ref{eq:linearmodel}) with $\bsbeta=(\beta_1,\beta_2,\beta_3)^\tran$ and $\bsX=(X_1,X_2,X_3)^\tran$ being a Gaussian random vector $\bsX \sim \dnorm(\bsmu,\bsSigma)$ with $\bsmu=(0,0,0)^\tran$. 
We assume that $X_1$ is independent from both $X_2$ and $X_3$, and that $X_2$ and $X_3$ may be correlated.
The covariance matrix reads:
$$\bsSigma=\left( \begin{array}{ccc}
\sigma_1^2 & 0 & 0 \\
0 & \sigma_2^2 & \rho \sigma_2 \sigma_3 \\
0 & \rho \sigma_2 \sigma_3 & \sigma_3^2 \end{array} \right)\, , \; -1 \leq \rho \leq 1\, .$$
We obtained the following analytical results.
\begin{equation}
\begin{array}{ll}
\sigma^2 & =  \displaystyle \mbox{Var} [f(\bsX) ] = \sum_{j=1}^3 \beta_j^2\sigma_j^2 + 2 \, \rho \, \beta_2 \beta_3 \sigma_2 \sigma_3 \, ,\\
Sh_1 & = \displaystyle  (\beta_1^2\sigma_1^2)/\sigma^2\, ,\\
Sh_2 & = \displaystyle  [\beta_2^2\sigma_2^2+\rho \, \beta_2\beta_3\sigma_2\sigma_3 + \frac{\rho^2}{2} (\beta_3^2\sigma_3^2-\beta_2^2\sigma_2^2)]/\sigma^2 \, ,\\
Sh_3 & = \displaystyle  [\beta_3^2\sigma_3^2+\rho \, \beta_2\beta_3\sigma_2\sigma_3 + \frac{\rho^2}{2} (\beta_2^2\sigma_2^2-\beta_3^2\sigma_3^2)]/\sigma^2 \, ,
\end{array}
\end{equation}
As expected, we have $\sum_{j=1}^3 Sh_j = 1$ and we see in $Sh_2$ and $Sh_3$ how the correlation effect is distributed in each index. 
In the case of fully correlated variables (i.e. $\rho=\pm 1$), we obtain $\displaystyle Sh_2=Sh_3=\left(\beta_2^2\sigma_2^2+\beta_3^2\sigma_3^2 +2 \rho \, \beta_2\beta_3\sigma_2\sigma_3\right)/\left(2\sigma^2\right)$.

\begin{figure}[!ht]
  \centering
	\includegraphics[width=0.49\textwidth]{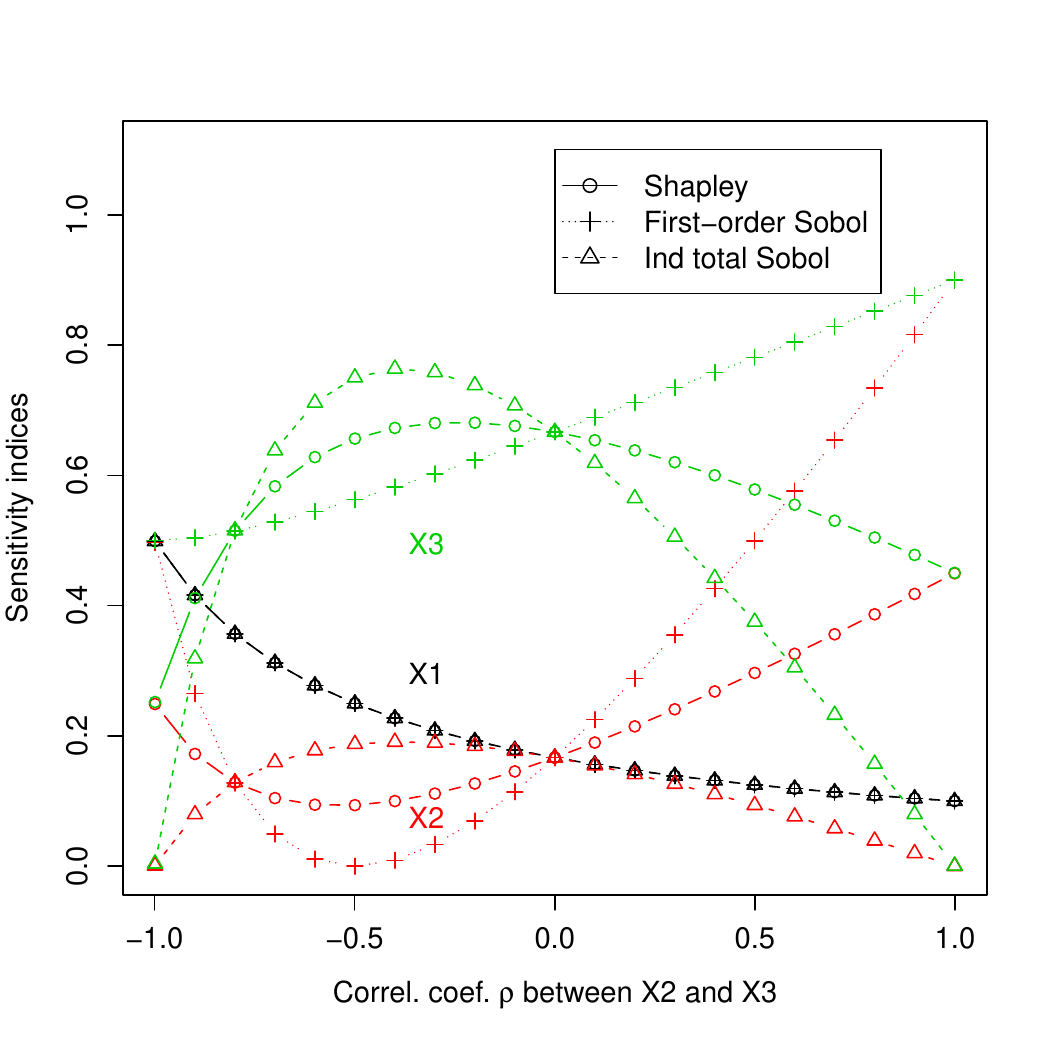}
		
\vspace{-0.4cm}
\caption{Sensitivity indices on the linear model with three Gaussian inputs.}
  \label{fig:modlin3d}
\end{figure}

We study the particular case $\beta_1=\beta_2=\beta_3=1$, $\sigma_1=\sigma_2=1$ and $\sigma_3=2$, for which in \cite{kuctar12}, the authors provide the formulas of full first-order and independent total Sobol' indices.
The analytical indices are depicted on Figure \ref{fig:modlin3d} as a function of the correlation coefficient $\rho$.
The Shapley effects are equal to the Sobol' indices in the absence of correlation, and then lie between the associated full first-order and independent total indices in the presence of correlation.
The ``sandwich effect'' is respected.
The effect of an increasing correlation (in absolute value) can be interpreted as an attractive effect both for the full Sobol' indices (here the first-order one) and the Shapley effect. 
However, for the Shapley effects, the contribution of the correlation is shared with each correlated variable.
This leads to the increase of one Shapley effect and the decrease of the other.
The Shapley effects allow an easy understanding of the influential inputs even when the Sobol' indices are not (when $ S_i> S_ {T_i} $).
As before, we see that the Shapley effects of two perfectly correlated variables are equal.

\subsection{Additive model with an interaction and three Gaussian inputs}\label{sec:lingausinterac}

In the Sections 3.2-3.4, we derived analytical results for which the Shapley effects were bounded by the full first-order and independent total indices. 
In the present section, we show that it is not always the case.
Let us define the model
\begin{equation}
Y=X_1+X_2X_3 
\end{equation} 
with 
\begin{equation*}
\bsX=\begin{pmatrix} 
X_1\\
X_2\\
X_3
\end{pmatrix}\sim
\dnorm_3\left(
\begin{pmatrix} 
0\\
0\\
0
\end{pmatrix}, \bsSigma\right)
\mbox{ and } \bsSigma=
\begin{pmatrix}
\sigma_1^2 & 0 & \rho \sigma_1 \sigma_3\\
0 & \sigma_2^2 & 0\\
\rho \sigma_1 \sigma_3 & 0 & \sigma_3^2
\end{pmatrix}, \; -1 \leq \rho \leq 1 \, .
\end{equation*}
It can be proven that $\sigma^2 S_1= \sigma_1^2$ and $\sigma^2 S_{T_1}=\sigma^2-\left(\sigma_2^2 \sigma_3^2 + \rho^2 \sigma_1^2 \right)=(1-\rho^2)\sigma_1^2$.
Recall that 
$$\sigma^2 Sh_1=\frac{1}{3}\left(\underline{\tau}^2_1-\underline{\tau}^2_{\emptyset}+\frac{1}{2} \left( \underline{\tau}^2_{12}-\underline{\tau}^2_2+\underline{\tau}^2_{13}-\underline{\tau}^2_3\right)+\sigma^2-\underline{\tau}^2_{23}\right) \, .$$
We thus get
\begin{equation}\label{shapcex}
\sigma^2 Sh_1=\sigma_1^2 (1-\frac{\rho^2}{2})+ \frac{\sigma_2^2\sigma_3^2}{6} \rho^2 \, .
\end{equation}
A straightforward computation yields $$S_{T_1} \leq Sh_1 \leq  S_1 \, .$$

We also get $S_2=0$, $\sigma^2 \, S_{T_2}=\sigma_2^2 \sigma_3^2$ and $\sigma^2 Sh_2=\frac{\sigma_2^2\sigma_3^2}{6}(3+ \rho^2)$. Thus
$$ S_2 \leq Sh_2 \leq  S_{T_2}\, .$$

Concerning the third input variable $X_3$, one gets $\sigma^2 S_3=\rho^2 \sigma_1^2$, $\sigma^2 S_{T_3}=(1- \rho^2) \sigma_2^2 \sigma_3^2$ and $\sigma^2 Sh_3= \frac{\rho^2 \sigma_1^2}{2} + \frac{\sigma_2^2\sigma_3^2}{6}(3-2 \rho^2)$.
Thus the two following assertions are equivalent:
$$S_3 \leq Sh_3 \leq S_{T_3} \,,$$
$$\rho^2 \sigma_1^2 \leq \frac{\sigma_2^2\sigma_3^2}{3}(3-4 \rho^2)\, .$$
The two following assertions are also equivalent:
$$ S_{T_3} \leq \frac{\phi_3}{\sigma^2} \leq S_3 \,,$$
$$\rho^2 \sigma_1^2 \geq \frac{\sigma_2^2\sigma_3^2}{3}(3-2 \rho^2)\, .$$
It also happens that $Sh_3$ is not comprised between $S_3$ and $S_{T_3}$.
The two following assertions are equivalent:
$$Sh_3 \geq \max \left( S_3, S_{T_3}\right) \,,$$
$$\frac37 \leq \rho^2 \leq \frac35 \,.$$

\begin{figure}[!ht]
  \centering
	\includegraphics[width=0.49\textwidth]{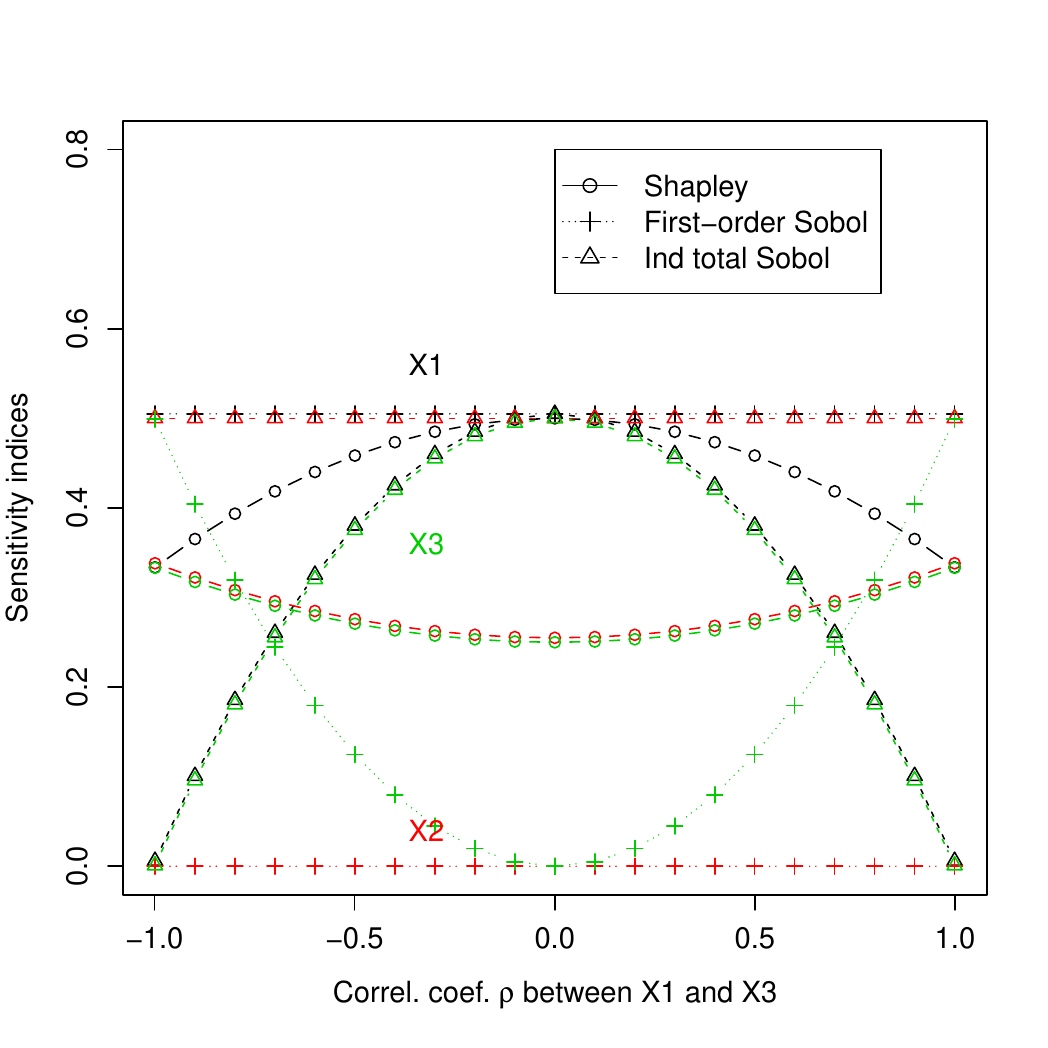}
		
\vspace{-0.4cm}
\caption{Sensitivity indices on the linear model with three Gaussian inputs and an interaction between $X_2$ and $X_3$.}
  \label{fig:modlinint3d}
\end{figure}

Figure \ref{fig:modlinint3d} illustrates the previous findings about the Sobol' indices and the Shapley effects for this model.
As expected, when the correlation coefficient $\rho$ belongs to two intervals, $[-0.775;-0.655]$ and $[0.655;0.775]$, the Shapley effects of $X_3$ are larger than the full first-order Sobol' indices and the independent total Sobol' indices.

\subsection{Three dimensional model with a block-additive structure}\label{sec:ishig}

We consider the following model: 
\begin{equation}\label{eq:addind}
Y=g(X_1,X_2)+h(X_3)\, ,
\end{equation}
which is called a ``block-additive'' structure.
We consider the general case where the vector $(X_1,X_2,X_3)^\tran$ is not restricted to a Gaussian vector. We only assume that the three inputs have finite variances and that $X_3$ is independent from $(X_1,X_2)$.
From the independence properties one has:
\begin{equation}\label{rem}
\sigma^2=\underline{\tau}^2_{12}+\underline{\tau}^2_3 \, , \; \;\;\underline{\tau}^2_{13}= \underline{\tau}^2_1 + \underline{\tau}^2_3  \;\mbox{ and }\; \underline{\tau}^2_{23}=\underline{\tau}^2_2+\underline{\tau}^2_3\, .
\end{equation}
From Equations \eqref{indclosed}, \eqref{indtotal2}, \eqref{indshapley} and (\ref{rem}) we get:
$$S_3=Sh_3=S_{T_3} \, .$$
We also get that, for $j=1,2$, the three following assertions are equivalent
$$S_j \leq Sh_j \,,$$
$$Sh_j \leq S_{T_j} \,,$$
$$\displaystyle \frac{\underline{\tau}^2_1+\underline{\tau}^2_2}{2} \leq \frac{\underline{\tau}^2_{12}}{2} \, .$$

We now consider, as in \cite{kuctar12}, the Ishigami function, a non-linear model involving interaction effects which writes:
\begin{equation}\label{eq:ishigami}
f(\bsX) = \sin(X_1) + 7 \sin(X_2)^2 +0.1 X_3^4 \sin(X_1)
\end{equation}
where $X_i \sim \dustd[-\pi,\pi] \; \forall i=1,2,3$ with a non-zero correlation $\rho$ between a pair of variables. 

\begin{figure}[!ht]
  \centering
	\includegraphics[width=0.49\textwidth]{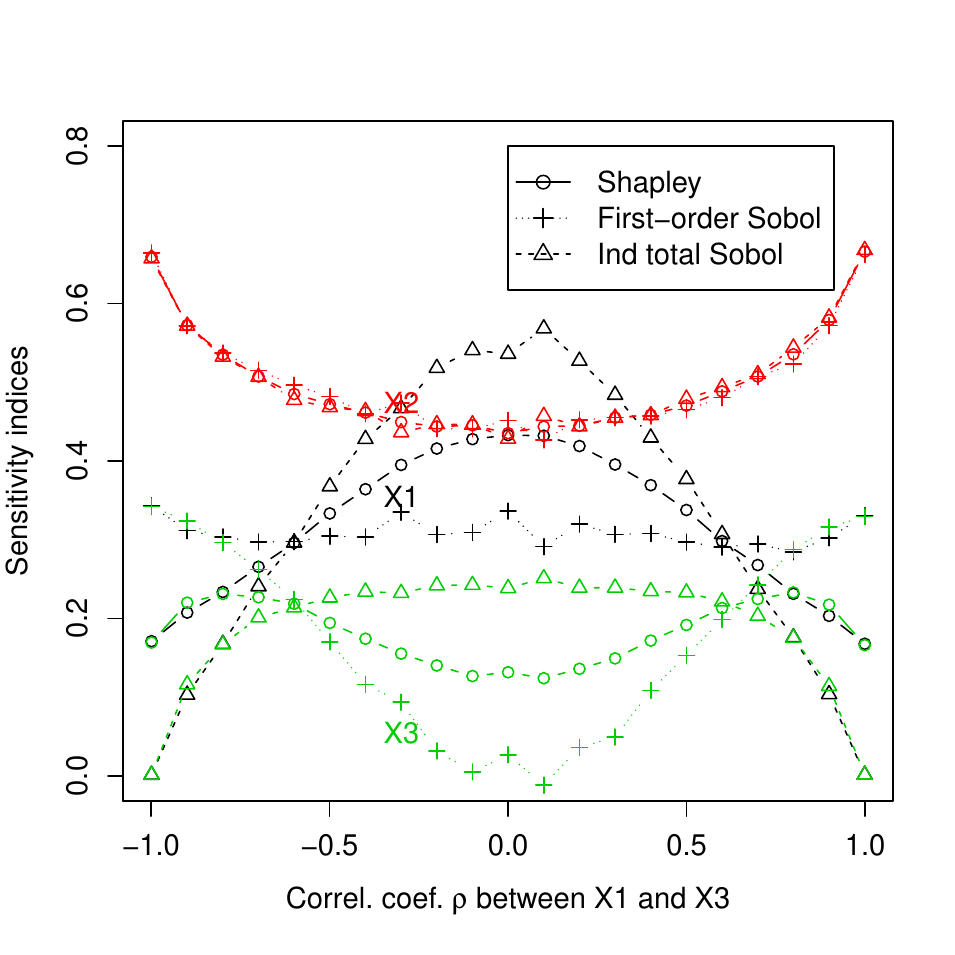}
	\includegraphics[width=0.49\textwidth]{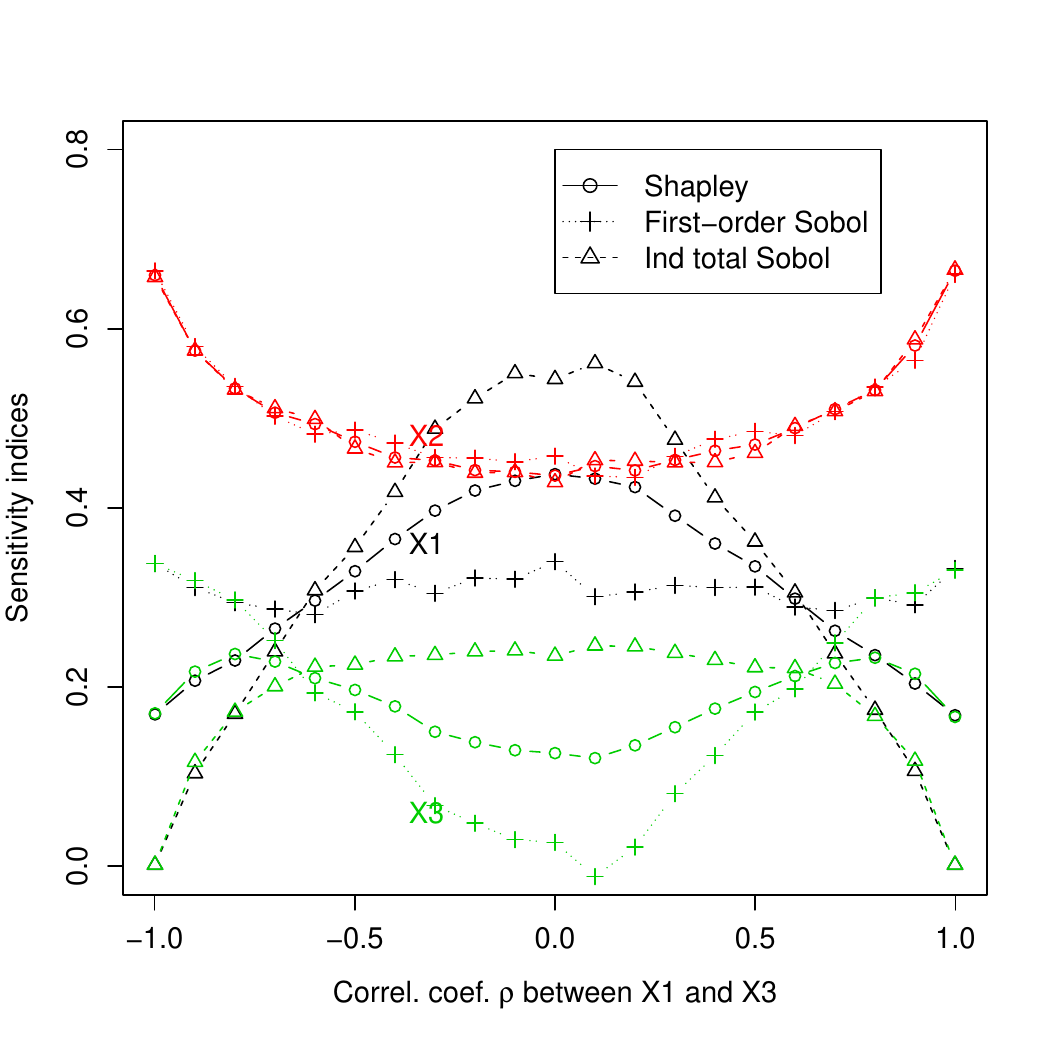}
		
\vspace{-0.4cm}
\centerline{(a) \hspace{8cm} (b)}

\caption{Sensitivity indices on the Ishigami function. (a) Exact permutation method with $N_o=2 \times 10^4$, $N_i=3$, $N_v=10^4$. (b) Random permutation method with $m=1.2\times 10^5$, $N_o=1$, $N_i=3$, $N_v=10^4$.}
  \label{fig:ishig13}
\end{figure}

Our study considers correlations between $X_1$ and $X_3$ only, $X_2$ being independent of $X_1$ and $X_3$. This model has a block-additive structure (as in Eq. \eqref{eq:addind} up to a permutation between $X_2$ and $X_3$).
{\color{black} The sensitivity measures depicted in Figure \ref{fig:ishig13} were obtained with the two numerical procedures explained in \ref{sec:direct}.
The computational cost ($C=7.3\times 10^5$) is the same between Fig. \ref{fig:ishig13} (a) and Fig. \ref{fig:ishig13} (b).
These two results have been provided in order to enhance our confidence on the correctness of the estimates as analytical values are unavailable.}
We observe that the ``sandwich effect'' is respected for $X_1$ and $X_3$. As $X_2$ is independent from the group $(X_1,X_3)$ and it has no interaction with that group, the Shapley index of $X_2$ equals both its full first-order and independent total indices.
These results confirm the general results discussed above for such a block-additive structure.
Moreover, the Shapley effects of $X_1$ and $X_3$ get closer as the correlation between them increases.

%
%
	%

\section{Practical issues and examples}\label{sec:estimation}

\subsection{Metamodel-based estimation}\label{sec:MM}

In this section, we consider a relatively common case in industrial applications where the numerical code is expensive in computational time.
As a consequence, it cannot be evaluated intensively (e.g. only several hundreds calculations are possible).
It is therefore not possible to estimate the sensitivity indices with direct use of the model.
Indeed, the Monte Carlo estimates of Sobol' indices require for each input several hundreds or thousands of model evaluations \cite{ioosal17,pritar17}. 
For the Shapley effects, an additional loop is required which increases the computational burden.

In this case, it is recommended to use a metamodel instead of the original numerical model in the estimation procedure.
A metamodel is an approximation of the numerical model, built on a learning dataset \cite{fanli06}.
The appeal to a metamodel is a usual engineering practice for estimating sensitivity indices \cite{legmar17}.
In the present work, we use Gaussian process-based metamodels (also called kriging) \cite{sacwel89,sanwil03} which have demonstrated in many practical situations having good predictive capacities (see \cite{marioo08} for example).
A Gaussian process model is defined as follows: 
\begin{equation} 
Y(\bsX) = h(\bsX) + Z(\bsX) , 
\end{equation}
where $h(\cdot)$ is a deterministic trend function (typically a multiple linear model) and $Z(\cdot)$ is a centered Gaussian process. 
The practical implementation details of kriging can be found in \cite{rougin12}.
We make the assumption that $Z$ is second-order stationary with a Mat\'ern $5/2$ covariance parameterized by the vector of its correlation lengths $\bstheta \in \real^d$ and variance $\sigma^2$. 
The hyperparameters $\sigma^2$ and $\bstheta$ are classically estimated by the maximum likelihood method on a learning sample comprising input/output of a limited number of simulations.
Kriging provides an estimator of $Y(\bsX)$ which is called the kriging predictor denoted by $\widehat{Y}(\bsX)$.
To quantify the predictive capability of the metamodel and to validate the predictor, the metamodel predictivity coefficient $Q^2$ is estimated by cross-validation or on a test sample \citep{marioo08}.
More precisely, the Gaussian process model gives the following predictive distribution:
\begin{equation}\label{eq:kriging}
\forall \bsX^\star \,, \quad \left(Y(\bsX^\star) \mid \bsy^N \right) \sim \dnorm \left( \widehat{Y}(\bsX^\star), \widehat{\sigma_Y^2}(\bsX^\star) \right) 
\end{equation}
where $\bsX^\star$ is a point of the input space not contained in the learning sample, $\bsy^N$ is the output vector of the learning sample of size $N$ and $\widehat{\sigma_Y^2}(\bsX)$ is the kriging variance that can also be explicitly estimated. 
In particular, the kriging variance $\widehat{\sigma_Y^2}(\bsX)$ quantifies the uncertainty induced by estimating $Y(\bsX)$ with $\widehat{Y}(\bsX)$. 

\begin{figure}[!ht]
  \centering
	\includegraphics[width=0.49\textwidth]{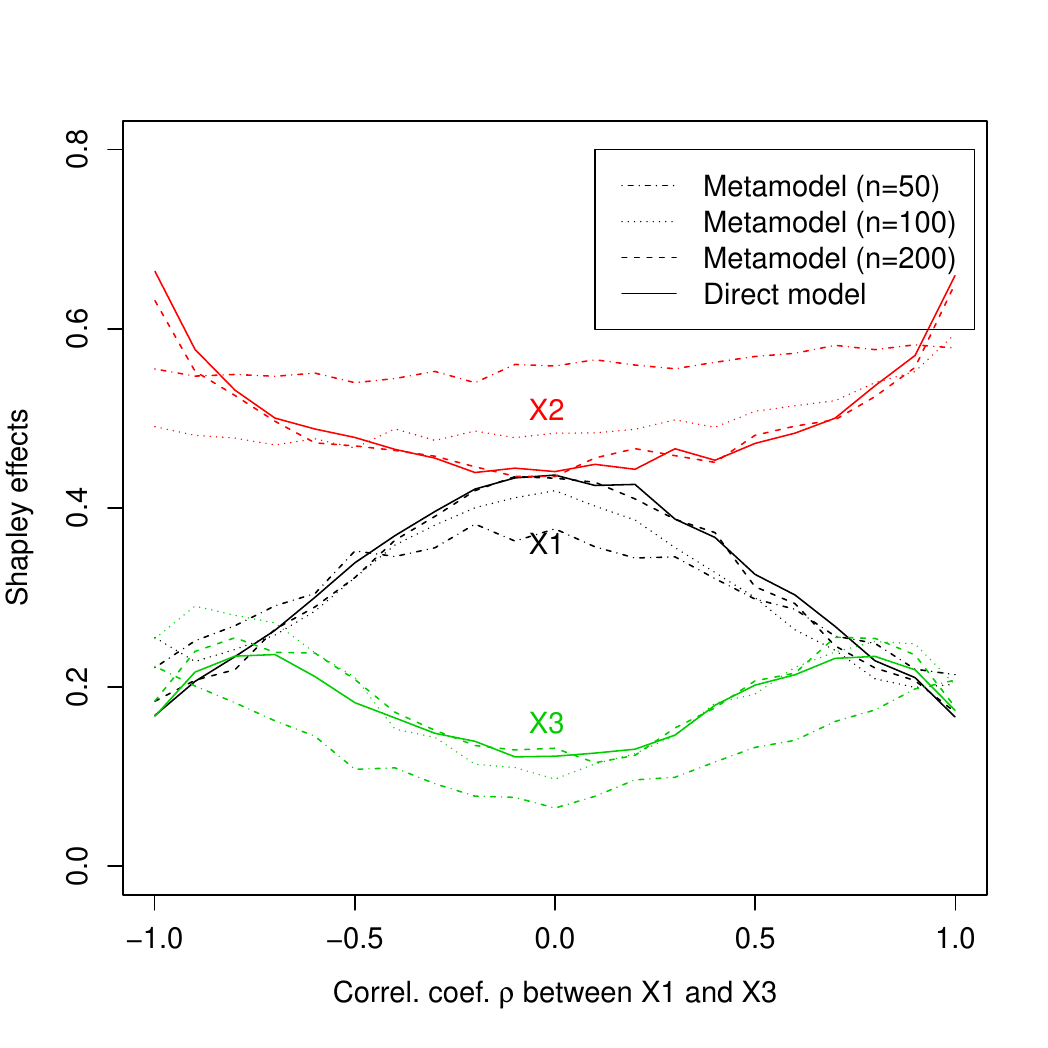}
	\includegraphics[width=0.49\textwidth]{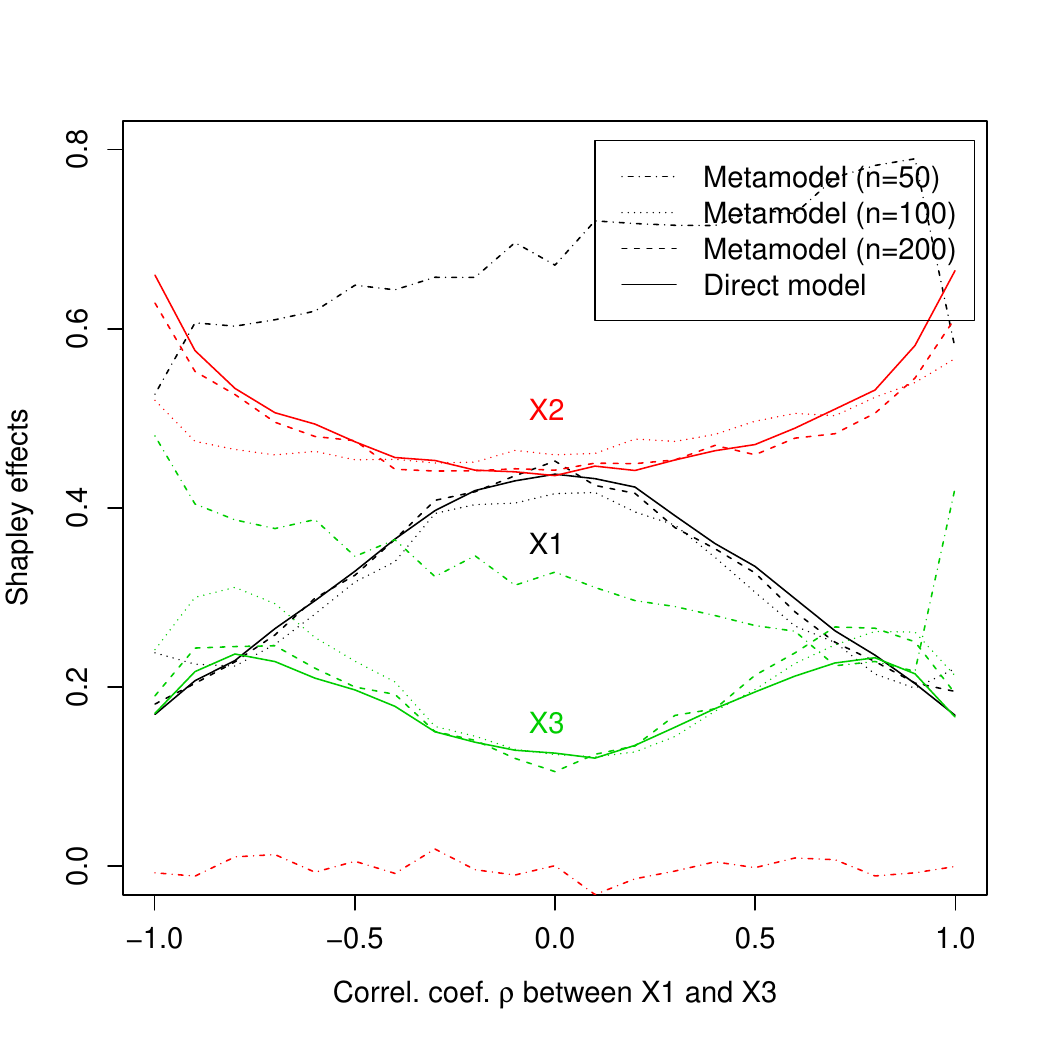}
	
\vspace{-0.4cm}
\centerline{(a) \hspace{8cm} (b)}

\caption{Results on the Ishigami function. Shapley effects estimated from $3$ different metamodels (corresponding to different sizes of the learning sample) and from the direct model (without metamodel). (a) Exact permutation method (with $N_o=2 \times 10^3$, $N_i=3$, $N_v=10^4$). (b) Random permutation method (with $m=1.2 \times 10^4$, $N_o=1$, $N_i=3$, $N_v=10^4$).}\label{fig:ishig}
\end{figure}

As an illustration, we study the Ishigami function (Eq. (\ref{eq:ishigami})) with a correlation coefficient $\rho$ between $ X_1 $ and $ X_3 $, on which \cite{kuctar12} studied the Sobol' indices (see also Section \ref{sec:ishig} of this paper).
When constructing the models, three different sizes $N$ of the learning sample ($50$, $100$ and $200$)  give three predictive coefficients for the kriging predictor: $ 0.78$, $0.88$ and $0.98$ respectively.
{\color{black}Figure \ref{fig:ishig} (a) and (b) were obtained with the two numerical procedures, explained in \ref{sec:direct}, applied on the kriging predictor.
The computational cost (total number of calls of the metamodel predictor) is the same in Fig. \ref{fig:ishig} (a) as in Fig. \ref{fig:ishig} (b)): $C=8.2\times 10^4$.
This high cost has no consequence as the metamodel predictor is evaluated instantly.
The results in Fig. \ref{fig:ishig} (a) and Fig. \ref{fig:ishig} (b) are similar.
They show that with a strong predictive metamodel ($Q^2=0.98$ with $N=200$), the estimations of the Shapley effects by the metamodel are satisfactorily accurate.
The precision of the estimated effects deteriorates rapidly with the decrease of the metamodel predictivity.
}

This example has just illustrated the need to have a sufficiently accurate metamodel in order to have precise estimates of Shapley effects.
{\color{black} However, replacing the original numerical model $f(\bsX)$ by a metamodel induces an additional error on the Shapley effect estimates due to the metamodel approximation.
As shown in \cite{legcan14} for Sobol' indices, it is possible to control this error thanks to the Gaussian process properties by using conditional Gaussian process simulations.
Such a development for Shapley effects is outside the scope of the present paper and has been recently proposed in \cite{beneli18}.
}

\subsection{Industrial application}\label{sec:example}

This application concerns a probabilistic analysis of an ultrasonic non-destructive control of a weld containing manufacturing defects.
Complex phenomena occur in such a heterogeneous medium during the ultrasonic wave propagation and a fine analysis to understand the effect of uncertain parameters is important.
The simulation of these phenomena is performed via the finite element code ATHENA2D, developed by EDF (Electricit\'e de France).
This code is dedicated to the simulation of elastic wave propagation in heterogeneous and anisotropic materials like welds.

\begin{figure}[!ht]
  \centering
	\includegraphics[width=0.59\textwidth]{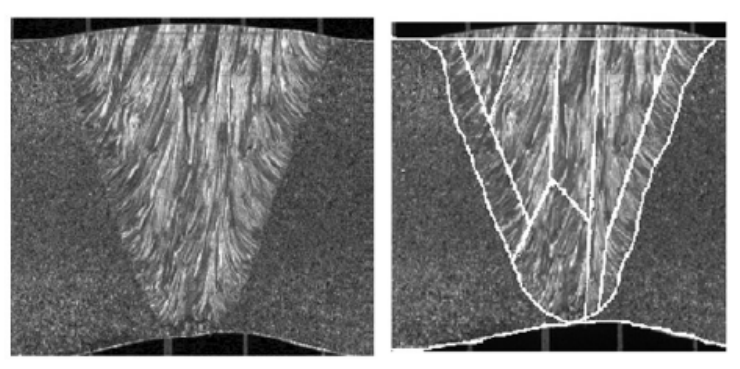}
	\includegraphics[width=0.4\textwidth]{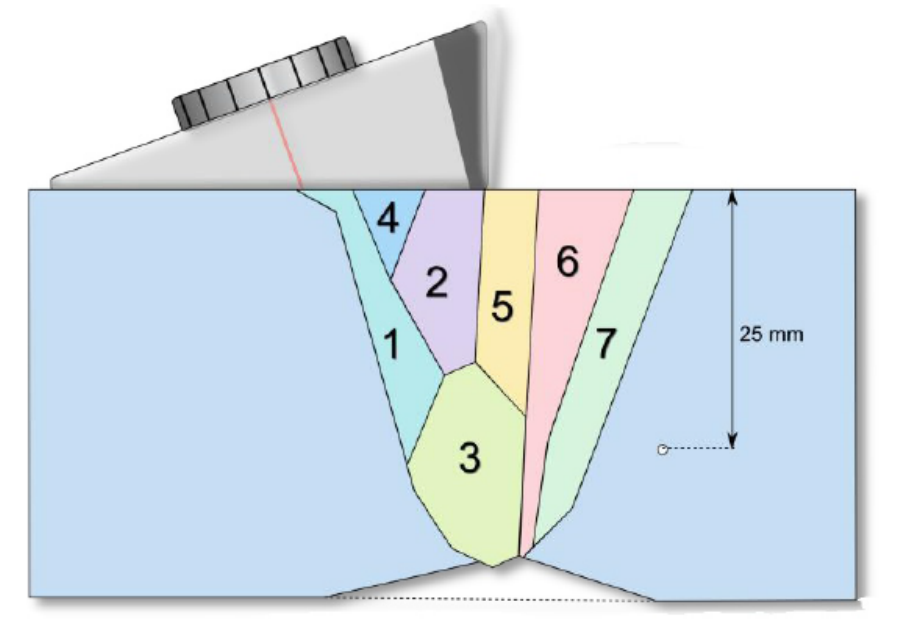}
	
\vspace{-0.4cm}
\caption{Metallographic picture (left), description of the weld in $7$ homogeneous domains (middle) and inspection configuration (right). From \cite{rupbla14}.}\label{fig:soudure}
\end{figure}

A first study \cite{rupbla14} has been realized with an inspection configuration aiming to detect a manufactured volumic defect located in a $40$ mm thick V grooveweld made of $316$L steel (Figure \ref{fig:soudure}).
The weld material reveals a heterogeneous and anisotropic structure. It is represented by a simplified model consisting of a partition of $7$ equivalent homogeneous regions with a specific grain orientation.
Eleven scalar input variables ($4$ elastic coefficients and $7$ orientations of the columnar grains of the weld inspections) have been considered as uncertain and modeled by independent random variables, each one associated to a probability density function.
The scalar output variable of the model is the amplitude of the defect echoes resulting from an ultrasonic inspection (maximum value on a so-called Bscan).
Uncertainty and sensitivity analysis (based on polynomial chaos expansion \cite{legmar17}) have then been applied from $6000$ Monte Carlo simulations of ATHENA2D in \cite{rupbla14}.
The sensitivity analysis has shown that almost all inputs are influential (only one input has a total Sobol' index smaller than $5\%$), that the interaction effects are non-negligible (approximately $30\%$ of the total variance) and that the orientations play a major role for explaining the amplitude variability.
The analysis confirms that an accurate determination of the micro-structure is essential in these simulation studies.
Finally, as a perspective of their work, the authors in \cite{rupbla14} explain that the real configuration has been strongly simplified by considering independent input random variables.
Indeed, due to the welding physical process, a dependence structure exists between the orientations, in particular between two neighboring domains (see Figure \ref{fig:soudure} right).

The purpose of the present study is then to perform a sensitivity analysis by using a more realistic probabilistic model for the input random variables.
Our SA setting is mainly a FF objective (see Section \ref{sec:SAsettings}): Which parameters can be fixed without impacting the predicted model response uncertainty?
Indeed, these SA results are expected to be useful in regards to the qualification process of the non-destructive control technique. 
As explained in Section \ref{sec:SAsettings}, the Shapley effects are well adapted to FF setting in the case of dependent inputs.

In our study, the probability distributions of all the inputs are considered Gaussian.
Their means and standard deviations are unchanged with respect to those of \cite{rupbla14}.
From physical models of welding process and solidification \cite{moyapf03}, engineers have been able to estimate the following correlation matrix between the $7$ orientations $(Or_1,\ldots,Or_7)$ of Figure \ref{fig:soudure} (right),
\begin{equation}\label{eq:Sigma}
\bsSigma=\left( \begin{array}{ccccccc}     
1    &0.80  &0.74  &0.69  &0.31  &0.23  &0.20 \\
0.80 &1     &0.64  &0.53  &0.59  &0.51  &0.46 \\
0.74 &0.64  &1     &0.25  &0.60  &0.57  &0.54 \\
0.69 &0.53  &0.25  &1     &-0.25 &-0.35 &-0.33 \\
0.31 &0.59  &0.60 &-0.25  &1     &0.96  &0.84 \\
0.23 &0.51  &0.57 &-0.35  &0.96  &1     &0.95 \\
0.20 &0.46  &0.54 &-0.33  &0.84  &0.95  &1        
	\end{array}\right)\,.
	\end{equation}

As only several hundreds of numerical simulations of ATHENA2D can be performed in the schedule time of the present study, our strategy consists of generating a space filling design in order to have a ``good'' learning sample for a metamodel building process. 
A Sobol' sequence of $N=500$ points has then been generated for the $d=11$ input variables on $[0,1]^d$.
After transformation of this sample to a sample of inputs which follow their physical scales and their joint probability density function, the corresponding $500$ runs of ATHENA2D have been computed.

\noindent {\bf Remark:} The $6000$ Monte Carlo simulations performed in the previous study \cite{rupbla14} were not stored, and thus could not be reused. As already mentioned, the metamodel built in that previous study was based on polynomial chaos expansion, and was not stored as well.

From the resulting $N$-size learning sample, a Gaussian process metamodel (parameterized as explained in Section \ref{sec:MM}) has then be fitted. 
We refer to \cite{legmar17} for a comparative study between metamodels based on polynomial chaos expansions and the one based on Gaussian processes.
We obtain a predictivity coefficient of $Q^2=87\%$.
This result is rather satisfactory, especially when it is compared to the predictivity coefficient obtained by a simple linear model ($Q^2=25\%$).
Moreover, the test on Ishigami function (Section \ref{sec:MM}) has shown that the estimation of Shapley effects with a metamodel of predictivity close to $90\%$ gives results rather close to the exact values. 

\begin{figure}[!ht]
  \centering
	\includegraphics[width=8cm, height=7cm]{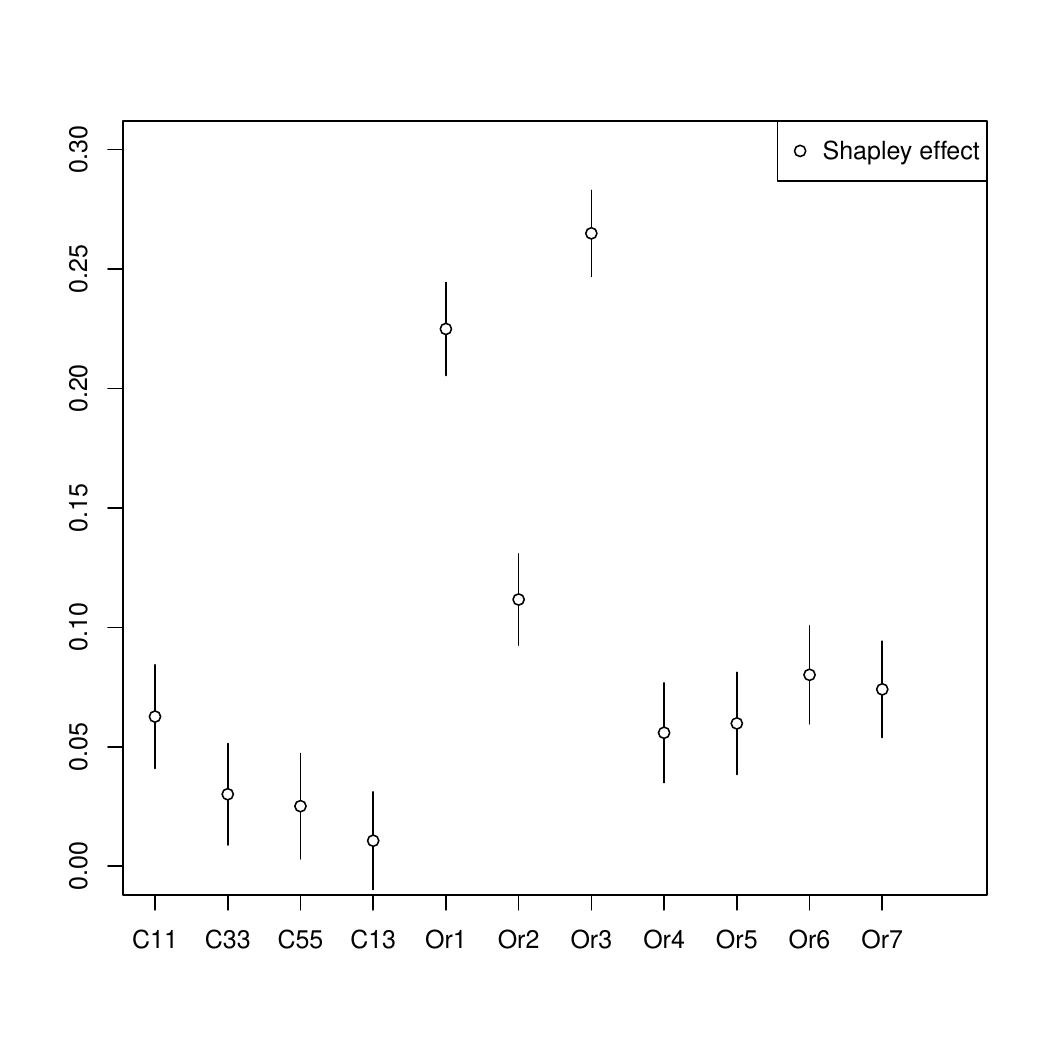}
	
\vspace{-0.4cm}
\caption{Shapley effects for the ultrasonic non-destructive control application.
The vertical bars represent the $95\%$-confidence intervals of each effect (uncertainty coming from the random permutation algorithm).}\label{fig:CNDUS}
\end{figure}

The Shapley effects are estimated by using the metamodel predictor instead of ATHENA2D (Section \ref{sec:MM}).
Due to the input dimension ($d=11$), the random permutation method is used with $m=10^4$, $N_i=3$, $N_o=1$ and $N_v=10^4$.
The cost is then $3\times {10}^5$ in terms of required metamodel evaluations.
It would be prohibitive with the ``true'' computer code ATHENA2D, but it is feasible by using the metamodel predictor.
Figure \ref{fig:CNDUS} gives the Shapley effects of the elasticity coefficients ($C_{11}$, $C_{13}$, $C_{33}$, $C_{55}$) and orientations ($Or_1$, $Or_2$, $Or_3$, $Or_4$, $Or_5$, $Or_6$, $Or_7$).
The lengths of the $95\%$-confidence interval (see Section \ref{sec:direct}) are approximately equal to $4\%$, which is sufficient to provide a reliable interpretation.
Note that the negative values of some Shapley effects are due to the central limit theorem approximation.

By visualizing the Shapley effects, we can propose a discrimination in four groups of inputs according to their degree of influence (note that such discrimination is questionable due to the residual uncertainties on Shapley effects):
\begin{itemize}
\item $Or_1$ and $Or_3$ whose effects are larger than $20\%$,
\item $Or_2$ whose effect is $11\%$,
\item $C_{11}$, $Or_4$, $Or_5$, $Or_6$ and $Or_7$ whose effects range between $6\%$ and $8\%$,
\item $C_{33}$, $C_{55}$ and $C_{13}$ whose effects are smaller than $3\%$.
The FF setting could be addressed with the inputs in this group.
\end{itemize}
To be convinced by this FF setting, the variance of the metamodel output when ($C_{33}$,$C_{55}$,$C_{13}$) are fixed is compared with the variance of the metamodel output when all the inputs vary.
{\color{black} These three inputs have been fixed to their mean values (note that they are independent of each other and independent of other inputs).}
The variance with all inputs is $3.774\times 10^{-22}$ and the variance with the three fixed inputs is $3.572\times 10^{-22}$.
As expected, the decrease of $5.3\%$ corresponds approximately to the sum of the Shapley effects of $C_{33}$, $C_{55}$ and $C_{13}$ (approximately $6\%$).


In the study of \cite{rupbla14} which did not take into account the correlation, $C_{33}$ and $Or_4$ have been identified as influential inputs (effects larger than $9\%$).
This result shows the importance of taking into account the dependence structure between inputs and the usefulness of the Shapley effects for FF setting in this case.
If we compare the (normalized) total Sobol' indices of \cite{rupbla14} and the Shapley effects of our study, taking into account the correlation has led to:
\begin{itemize}
\item an increase in sensitivity indices for $Or_1$, $Or_2$ and $Or_3$,  
\item a decrease in sensitivity indices for $Or_7$ , 
\item similar sensitivity indices for $Or_4$, $Or_5$ and $Or_6$.
\end{itemize}
By looking at the input correlation matrix (Eq. (\ref{eq:Sigma})), we remark that we can distinguish two groups of inputs as a function of their correlation degrees: $(Or_1,Or_2,Or_3,Or_4)$ and  $(Or_5,Or_6,Or_7)$.
We observe the homogeneity of the correlation structure effects: the inputs inside the first group correspond to an increase (or a stability) in sensitivity indices whereas the inputs inside the second group correspond to a decrease (or a stability) in sensitivity indices.

\begin{figure}[!ht]
  \centering
	\includegraphics[width=14cm, height=9cm]{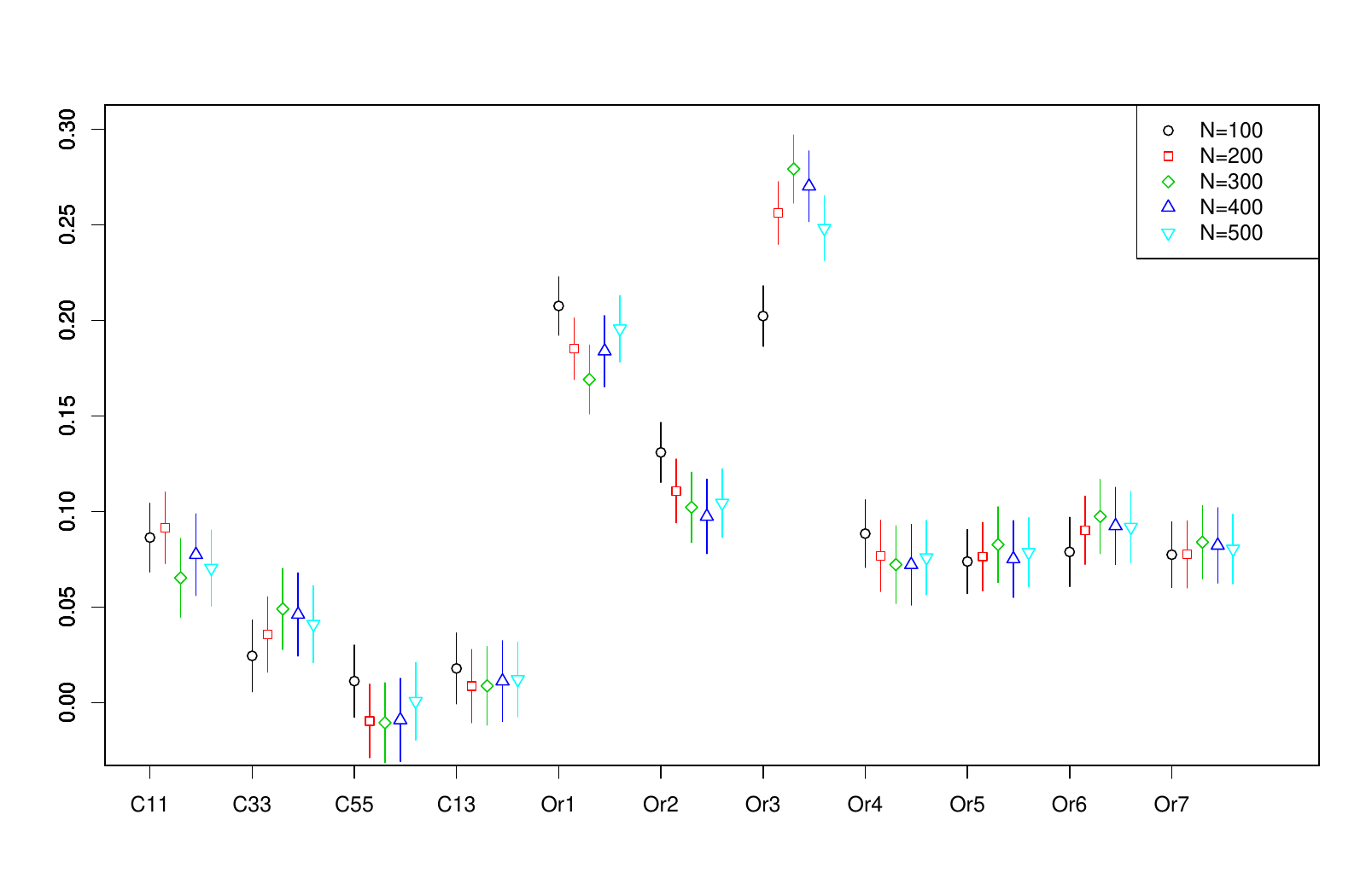}
	
\vspace{-0.4cm}
\caption{Convergence of the Shapley effects with respect to the metamodel learning sample size.
The vertical bars represent the $95\%$-confidence intervals of each effect (uncertainty coming from the random permutation algorithm).}\label{fig:CNDUS_cv}
\end{figure}

{\color{black} To confirm the results obtained in Figure \ref{fig:CNDUS}, a convergence study of the Shapley effect estimates with respect to the metamodel learning sample size is shown in Figure \ref{fig:CNDUS_cv}. 
The learning sample size varies from $N=100$ to $N=500$ which is the total available budget (corresponding to the results of Figure \ref{fig:CNDUS}).
The sample points are sequentially taken in the total sample. As the input design is a Sobol' sequence, it means that the input samples are nested within each other so as to fill the input space in a sequential and complementary manner.
Figure \ref{fig:CNDUS_cv} shows that, even for small sample sizes ($N=100$ and $N=200$) which yield coarse metamodels (negative estimated predictive coefficients $Q_2$), the estimates are rather close to the estimates obtained at $N=500$.
Then, for $N=300$ and $N=400$ (which gives $Q_2=0.29$ and $Q_2=0.55$), the convergence seems to have been reached for almost all the Shapley effect estimates.
We can conclude in a good confidence in the values given by the estimates at $N=500$ ($Q_2=0.87$), which could be slightly refined by increasing the learning sample size.
}

\section{Conclusion}

In many applications of global sensitivity analysis methods, it is common that the input variables have a known statistical dependence structure or that the input space is constrained to a non-rectangular region. 
In this paper we considered two answers to that issue: the Shapley effects (a normalized version of the variance-based Shapley values proposed in \cite{owe14} in the framework of sensitivity analysis) and the methodology developed in \cite{martar15}.
The latter suggests the joint analysis of full and independent first-order and total indices to analyze the sensitivity of a model to dependent inputs. 
In the present paper, we conducted a comparative analysis between Shapley effects on one side and full first-order and independent total indices on the other.
{\color{black}It is well known that full first-order indices do not contain enough information to allow factor fixing. It  also seems clear, from analytical solutions obtained with linear models and Gaussian variables, that even in the framework of independent inputs, total Sobol' indices encapsulate redundant information, while Shapley effects, due to the equitable principle on which they rely, provide normalized indices, which can be used for the FF Setting.} 
Comparisons of Shapley effects with the complementary independent first-order and full total indices are currently under investigation.

We have also illustrated the convergence of two numerical algorithms for estimating Shapley effects.
Our preliminary comparative study between Shapley effects and Sobol' indices is completed by the computations of $S_{T_{(i)}}$ and $S_{(i)}^{\textup{ind}}$ in \cite{beneli18} (see Section \ref{sob:ind}).
The studied algorithms depend on various parameters: $N_i$ (conditional variance estimation sample size), $N_o$ (expectation estimation sample size), $N_v$ (output variance estimation sample size) and $m$ (random permutation number). 
It would be interesting to investigate further the response of the algorithms to these different parameters and to derive empirical and asymptotic confidence intervals for the Shapley effects estimates. 
Introducing a sequential procedure in the random permutation algorithm, in order to increase $m$ until there is a sufficient precision in the Shapley effects, also seems  promising.
Moreover, it would be important in a future work to consider the estimation algorithm capabilities on more complex dependence structures than the pairwise cases exclusively discussed in the present paper.

Finally, we have shown the relevance of using a metamodel (here the Gaussian process predictor) in industrial situations where the computer model is too time consuming to be evaluated thousands of times for the previous algorithms to be applied. 
Future work (started in \cite{beneli18}) will consist of developing an algorithm exploiting the complete structure of the Gaussian process allowing to infer the error due to this approximation (see \cite{jannod14,jankle13,legcan14} for the Sobol' indices and \cite{delmar16} for the derivative-based global sensitivity measures).


\section*{Acknowledgments}

We are grateful to the useful and numerous comments of the two reviewers which greatly helped to improve the manuscript.
We thank G\'eraud Blatman who performed the computations on the ATHENA2D model, Roman Sueur for helpful discussions about Shapley effects, Chu Mai for his remarks, and Lynn Ferrieu for her help in proofreading the paper.
Numerical estimation of Shapley effects and Sobol' indices have been realized using the \textsf{sensitivity} package of the \textsf{R} software.
Thanks to Eunhye Song, Barry L. Nelson and Jeremy Staum for providing preliminary versions of the Shapley effects estimation functions, which are available at present in the ``sensitivity'' package of the R software.



\appendix

\section{Estimation of Shapley effects by direct sampling algorithms}\label{sec:direct}

The authors in \cite{sonnel16} propose two algorithms for estimating the Shapley effects from formula (\ref{indshapley}) with 
\begin{equation}\label{eq:cost}
c(u) = \frac{\mathbb{E}(\mbox{Var}[Y|\bsX_u])}{\mbox{Var}(Y)}
\end{equation}
being the cost function (which has been shown to be more efficient than the variance of the conditional expectation).
The first algorithm traverses all possible permutations between the inputs and is called the ``Exact permutation method''.
The second algorithm consists of randomly sampling some permutations of the inputs and is called the ``Random permutation method''. 
The latter is to be preferred when the overall permutation set is too large to be considered.

For each iteration of the inputs' permutations loop, a conditional variance expectation must be computed.
The cost $C$, in terms of model evaluations, of these algorithms is then the following \cite{sonnel16}:
\begin{enumerate}
\item Exact permutation method: $C = N_i N_o^{\textup{exa}} d! (d-1) + N_v$, with $N_i$ the inner loop size (conditional variance) in Eq. (\ref{eq:cost}), $N_o^{\textup{exa}}$ the outer loop size (expectation) in (\ref{eq:cost}) and $N_v$ the sample size for the variance computation (denominator in (\ref{eq:cost}));
\item Random permutation method: $C = N_i N_o^{\textup{rand}} m (d-1) + N_v$, with $m$ the number of random permutations for discretizing the principal sum in (\ref{indshapley}), $N_i$ the inner loop size, $N_o^{\textup{rand}}$ the outer loop size and $N_v$ the sample size for the variance computation.
\end{enumerate}
The $(d-1)$ terms that appear in the computational costs come from the fact that $(d-1)$ Shapley effects are estimated while the last Shapley effect is estimated by using the sum-to-one property.
Note that the full first-order Sobol' indices (Eq. (\ref{indclosed})) and the independent total Sobol' indices (Eq. (\ref{indtotal2})) are also estimated by applying these algorithms.

From theoretical arguments, the authors in \cite{sonnel16} have shown that the near-optimal values of the sizes of the different loops are the following:
\begin{itemize}
\item $N_i=3$ and $N_o^{\textup{exa}}$ as large as possible for the exact permutation method,
\item $N_i=3$, $N_o^{\textup{rand}}=1$ and $m$ as large as possible for the random permutation method.
\end{itemize}
With these choices, the value of $m$ which leads to the same numerical cost for the two algorithms is:
\begin{equation}\label{eq:costequival}
m = N_o^{\textup{exa}} d! \;.
\end{equation}
In practical applications with large $d$, this choice is not realistic and the random permutation method is applied with a much smaller value for $m$ (see Section \ref{sec:example}).

The exact permutation algorithm with fixed $N_i$ is consistent as $N_o^{\textup{exa}}$ tends to infinity. 
The random permutation one with fixed $N_i$ and $N_o^{\textup{rand}}$ is consistent as the number of sampled permutations, $m$, tends to infinity. Indeed, both algorithms are based on unbiased estimators of $\mbox{Var}(Y)$ and $Sh_i \times \mbox{Var}(Y)$, whose variance tends to zero (see Appendix A, Equations (18) and (22) in \cite{sonnel16} for more details). From Theorem 3 in \cite{sonnel16}, we know that the variance of the estimator of $Sh_i$ obtained from the random permutation algorithm is bounded by $\left(\mbox{Var}(Y)\right)^2/m$. More intuitively, it seems reasonable to think that the difficulty to estimate $Sh_i$ is related on the effective dimension of the model as well as on the complexity of the dependence structure of the inputs.

\begin{figure}[!ht]
  \centering
	\includegraphics[width=0.49\textwidth]{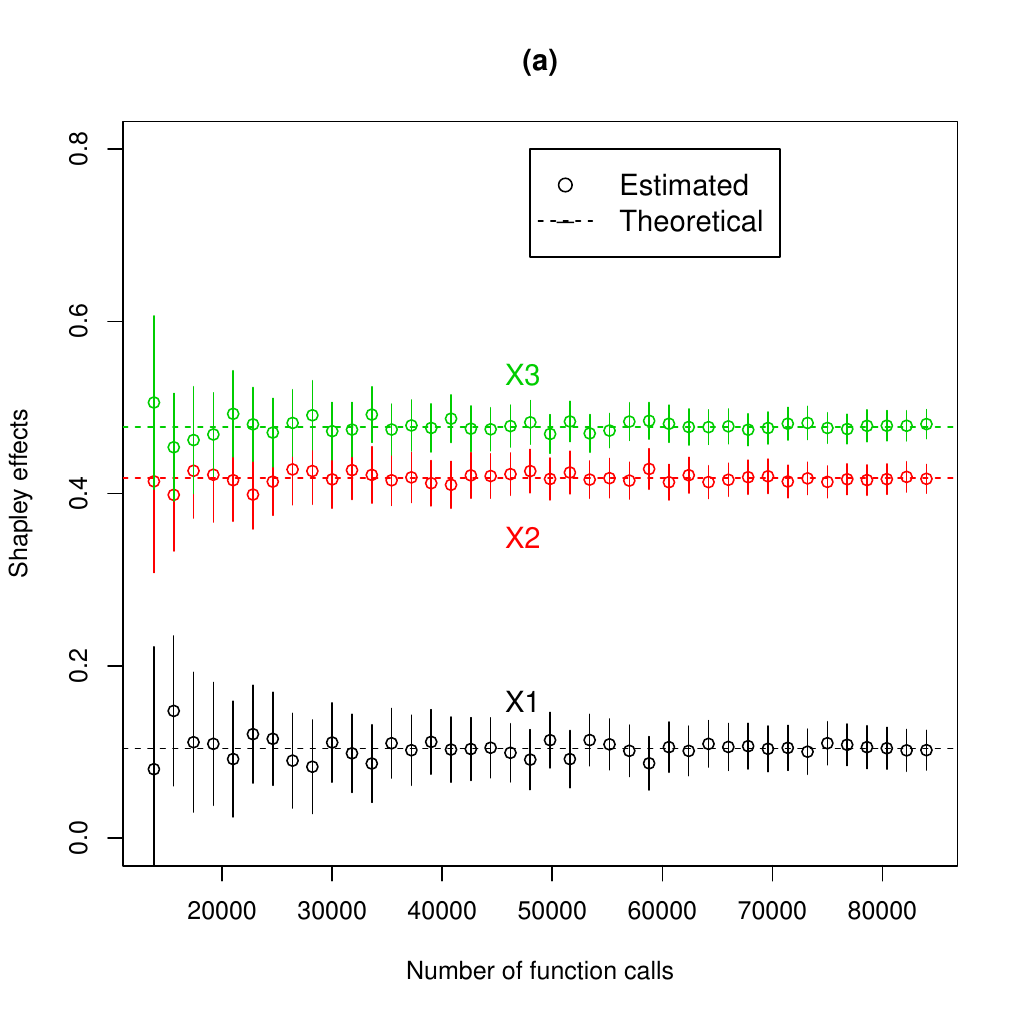}
	\includegraphics[width=0.49\textwidth]{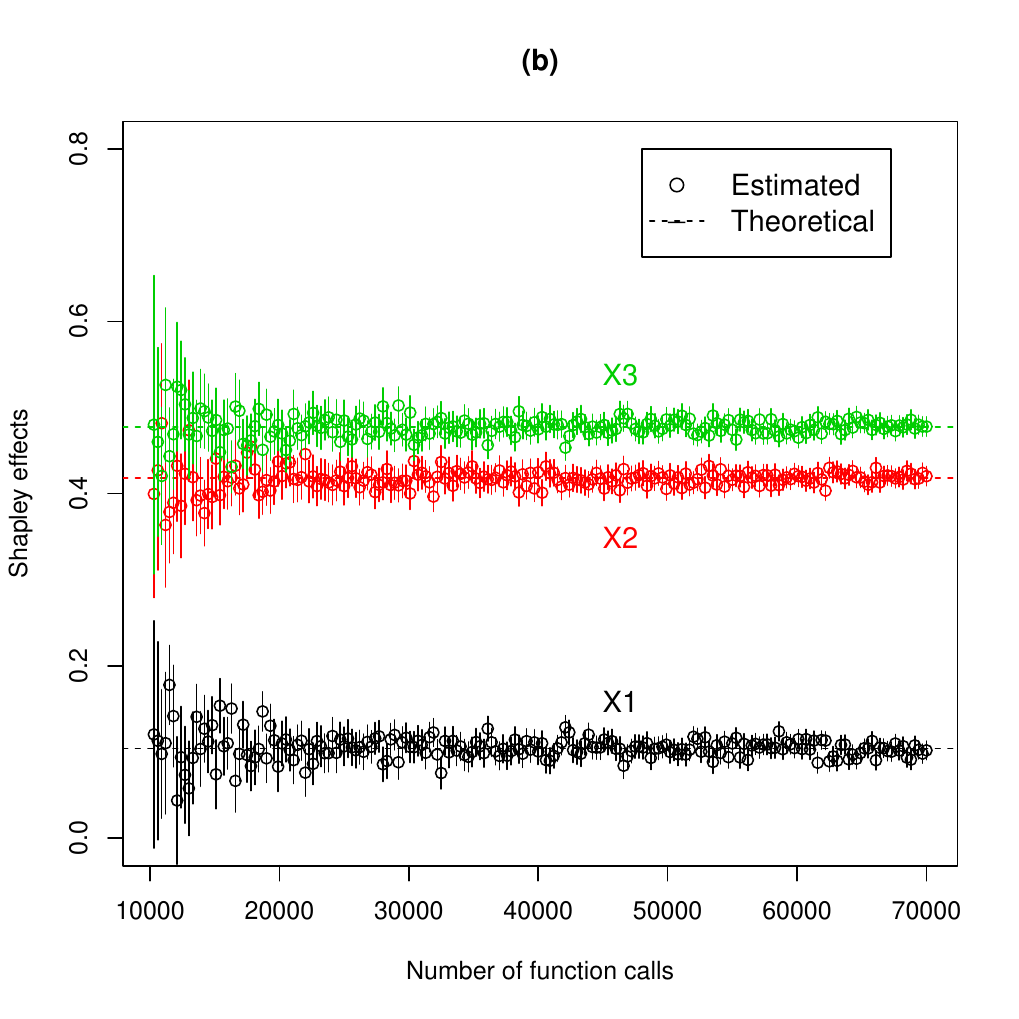}
	
\vspace{-0.4cm}
\caption{Numerical estimates by the exact permutation method (a) and the random permutation method (b) of the Shapley effects on the linear model with $3$ Gaussian inputs. }\label{fig:cvg_modlin}
\end{figure}

{\color{black} To illustrate these numerical estimators (with $N_i=3$ and $N_o^{\textup{rand}}=1$), we consider the linear model (Eq. (\ref{eq:linearmodel})) with $3$ Gaussian inputs of Section \ref{sec:lin3gau} with $\beta_1=\beta_2=\beta_3=1$, $\sigma_1=\sigma_2=1$, $\sigma_3=2$. 
We set $\rho=0.9$ (correlation between $X_2$ and $X_3$).
On Figure \ref{fig:cvg_modlin} (a) (resp. (b)), the results of the exact permutation method with $N_o^{\textup{exa}}$ ranging from $50$ to $2000$ (resp. random permutation method with $m$ ranging from $50$ to $10000$) are plotted versus the total cost $C$ (number of model evaluations). }
For the random permutation method, the error bars were obtained from the central limit theorem on the permutation loop (Monte Carlo sample of size $m$) and then by taking two times the standard deviation of the estimates ($95\%$ confidence intervals). Similarly, for the exact permutation method, the error bars were obtained from the central limit theorem on the outer loop (Monte Carlo sample of size $N_o^{\textup{exa}}$). 
In the both cases, we observe the convergence of the estimated values toward the exact values as $N_o^{\textup{exa}}$ (reps. $m$) increases.

{\color{black} Our goal in this appendix was just to illustrate the numerical convergence of the Shapley effect estimates.}
Further analysis about the choice  of $m$, $N_o^{\textup{exa}}$, $N_o^{\textup{rand}}$ and $N_i$ is nevertheless necessary.

%

\end{document}